\documentclass[11pt]{article}

\usepackage{color}
\usepackage{latexsym}
\usepackage{amssymb}
\usepackage{amsmath}
\usepackage{graphicx}
\usepackage{enumerate}
\usepackage{hyperref}

\newtheorem{Theorem}{Theorem}[part]
\newtheorem{Definition}{Definition}[part]
\newtheorem{Proposition}{Proposition}[part]

\newtheorem{Lemma}{Lemma}[part]
\newtheorem{Corollary}{Corollary}[part]
\newtheorem{Remark}{Remark}[part]

\renewcommand{\theDefinition}{\thesection.\arabic{Definition}}

\renewcommand{\theLemma}{\thesection.\arabic{Lemma}}

\renewcommand{\theequation}{\thesection.\arabic{equation}}

\def \Frac{\displaystyle\frac}

\def \N{\mathbb{N}}
\def \R{\mathbb{R}}

\def \E{\mathbb{E}}
\def \F{\mathbb{F}}

\def \P{\mathbb{P}}

\def \Ac{{\cal A}}

\def \Cc{{\cal C}}

\def \Fc{{\cal F}}
\def \Gc{{\cal G}}
\def \Hc{{\cal H}}

\def \Lc{{\cal L}}

\def \Oc{{\cal O}}
 \def \Nc{{\cal N}}
\def \Sc{{\cal S}}
\def \Tc{{\cal T}}
\def \Vc{{\cal V}}

\def \Vc{{\cal V}}

\def \eps{\varepsilon}

\def \ep{\hbox{ }\hfill$\Box$}

\def\Dt#1{\Frac{\partial #1}{\partial t}}
\def\Dth#1{\Frac{\partial #1}{\partial \theta}}

\def\Dp#1{\Frac{\partial #1}{\partial p}}

\def\Dpp#1{\Frac{\partial^2 #1}{\partial p^2}}

\def\Dpp#1{\Frac{\partial^2 #1}{ \partial p^2}}
\def\Dp#1{\Frac{\partial #1}{ \partial p}}

\def\Dpp#1{\Frac{\partial^2 #1}{\partial p^2}}

\def\reff#1{{\rm(\ref{#1})}}

\def\beqs{\begin{eqnarray*}}
\def\enqs{\end{eqnarray*}}
\def\beq{\begin{eqnarray}}
\def\enq{\end{eqnarray}}

\begin{document}

\title{Optimal portfolio liquidation  with execution cost and risk\thanks{We would like to thank Bruno Bouchard for useful comments. We also thank participants at the {\it Istanbul workshop on Mathematical Finance} in may 2009, for  relevant remarks.}}

\author{Idris KHARROUBI
             \\\small  Laboratoire de Probabilit\'es et
             \\\small  Mod\`eles Al\'eatoires
             \\\small  CNRS, UMR 7599
             \\\small  Universit\'e Paris 7, 
             \\\small  and CREST, 
             \\\small  e-mail: kharroubi@ensae.fr
             \and
             Huy\^en PHAM
             \\\small  Laboratoire de Probabilit\'es et
             \\\small  Mod\`eles Al\'eatoires
             \\\small  CNRS, UMR 7599
             \\\small  Universit\'e Paris 7, 
             \\\small  CREST, and 
             \\\small   Institut Universitaire de France
              \\\small  e-mail: pham@math.jussieu.fr
             }

\maketitle

\begin{abstract}
We study the optimal portfolio liquidation problem over a finite horizon in a limit order book with bid-ask spread and 
temporary market price impact penalizing speedy execution trades. We use a continuous-time modeling framework, but 
in contrast with previous related papers (see e.g.  \cite{rogsin08} and   \cite{schsch07}),  
we do not assume continuous-time trading strategies. We consider instead real trading 
that occur  in discrete-time, and  this is formulated as an impulse control problem under a solvency constraint, including the lag variable tracking the time interval between trades.  A first  important result of our paper is to show that nearly optimal execution strategies in this context lead actually to a finite number of  trading  times, and this holds true without assuming ad hoc any  fixed transaction fee.   
Next, we derive the dynamic programming quasi-variational 
inequality  satisfied by the value function in the sense of constrained viscosity solutions.  We also introduce a family of value functions converging to  our value function, and which is characterized 
as the unique constrained viscosity solutions of an approximation of our dynamic programming equation.  This convergence result  is useful  
for numerical purpose, postponed in a further study.  
\end{abstract}

\vspace{5mm}


\vspace{3mm}

\noindent {\bf MSC Classification (2000) : }93E20, 91B28, 60H30, 49L25.

\newpage

\section{Introduction}

Understanding  trade execution strategies  is a key  issue for financial market practitio\-ners, and has attracted a growing attention from the academic researchers.  An important  pro\-blem  faced by stock traders is  how to liquidate large block orders of shares.   This is a challenge due to the following  dilemma.  By trading quickly,  the investor is subject to higher costs due to market impact reflecting the depth of the limit order book.  
Thus, to minimize price impact, it is generally beneficial to break up a large order into smaller blocks.  However, more gradual trading over time results in higher risks since the asset value   can vary  more  during the investment horizon in an uncertain environment.  There has been recently a considerable interest in the literature on such  liquidity effects, taking into account permanent and/or  temporary price impact, 
and problems of this type were studied by  Bertsimas and Lo \cite{berlo98},  Almgren and Criss \cite{almcri01},  Bank and Baum \cite{banbau04}, Cetin, Jarrow and Protter \cite{cetjarpro04},  Obizhaeva and Wang \cite{obiwan05},  He and Mamayski \cite{hemam04},  
Schied an Sch\"oneborn \cite{schsch07},  Ly Vath, Mnif and Pham \cite{lyvmnipha07}, Rogers and Singh \cite{rogsin08}, and  Cetin, Soner and Touzi \cite{cetsontou09}, to mention some of them.

There are essentially two popular formulation types  for  the optimal trading problem in the literature: discrete-time versus continuous-time.  
In the discrete-time formulation,  we may distinguish papers considering that trading take place at fixed deterministic times (see \cite{berlo98}),  at exogenous random discrete times given for example by the jumps of a Poisson process 
(see \cite{phatan08}, \cite{baylud09}), or at  discrete times decided optimally by the investor through an impulse control formulation 
(see \cite{hemam04} and \cite{lyvmnipha07}).  In this last case, one usually assumes the existence of  a fixed transaction cost paid at each trading 
in order to ensure that strategies do not accumulate in time and  occur really  at discrete points in time 
(see e.g. \cite{kor98} or  \cite{okssul02}).  The continuous-time trading formulation is not realistic in practice, 
but is commonly used  (as in \cite{cetjarpro04}, \cite{schsch07} or \cite{rogsin08}), 
due to the tractability and powerful theory of the stochastic calculus  typically illustrated by It\^o's formula.  In a perfectly liquid market without transaction cost and market impact, continuous-time trading is often justified by arguing that it is a  limit approximation of discrete-time trading when the time step goes to zero. However, one may question  the validity of such assertion  in the presence of liquidity effects.

In this paper, we propose a  continuous-time framework taking into account the main liquidity features and risk/cost tradeoff of portfolio execution:  
there is a bid-ask spread in the limit order book, and temporary market price impact penalizing rapid execution trades. However, in contrast with 
previous related papers (\cite{schsch07} or \cite{rogsin08}), we do not assume continuous-time trading strategies. We consider instead real trading that take place in discrete-time, and without assuming ad hoc any fixed transaction cost,  in accordance with the practitioner literature.  
Moreover, a key issue  in line of the banking regulation and solvency constraints  is to define in an economically meaningful way  
the portfolio value of a position in stock at any time, and this is addressed in our modelling.  These issues are  formulated conveniently through an impulse control problem  including the lag variable tracking the time interval between trades. Thus, we combine the advantages of the stochastic calculus techniques, and the realistic modeling of portfolio liquidation.  In this context, we study the optimal 
portfolio liquidation problem over a finite horizon: the investor seeks to unwind an initial position in stock shares by maximizing 
his expected utility from terminal liquidation wealth, and under a natural economic solvency constraint involving the liquidation value of a portfolio.

A first important result  of our paper is to show that that nearly optimal execution strategies in this modeling lead actually to a finite number of trading times.  While most models dealing with trading strategies via an impulse control formulation assumed fixed transaction cost in order to justify a posteriori the discrete-nature of trading times, we prove here that discrete-time trading appear naturally as a consequence of 
liquidity features represented by  temporary price impact and bid-ask spread.  
Next, we derive the dynamic programming quasi-variational inequality (QVI) satisfied by the value function in the sense of 
constrained viscosity solutions in order to handle state constraints.  There are some technical difficulties related to the nonlinearity of the impulse transaction function induced by the market price impact, and the non smoothness of the solvency boundary. 
In particular, since we do not assume a fixed transaction fee, which precludes the existence of a strict supersolution to the QVI,  
we can not prove directly a comparison principle (hence a uniqueness result) for the  QVI.  We then consider  two types of  approximations by introducing families of value functions converging to our ori\-ginal value function, and which are characterized as unique constrained viscosity solutions to their dynamic programming equations. This convergence result is useful for numerical purpose, postponed in a further study.

The plan of the paper is organized as follows. Section 2 presents the details of the model and formulates the liquidation problem. In Section 3, we 
show some interesting economical and mathematical properties of the model, in particular the finiteness of the number of trading strategies under illiquidity costs. Section 4 is devoted to the dynamic programming and viscosity properties of the value function to our impulse control problem. 
We propose in Section 5 an approximation of the original problem by considering small fixed tran\-saction fee. Finally, Section 6 describes another approximation of the model with utility penalization by small cost. As a consequence, we obtain that our initial  value function   is characterized as the  minimal constrained viscosity solution to its dynamic programming QVI.

\section{The  model and liquidation problem}

\setcounter{equation}{0} \setcounter{Assumption}{0}
\setcounter{Theorem}{0} \setcounter{Proposition}{0}
\setcounter{Corollary}{0} \setcounter{Lemma}{0}
\setcounter{Definition}{0} \setcounter{Remark}{0}
 \setcounter{figure}{0}

We consider a financial market  where an  investor has to  
liquidate an initial  position of $y$ $>$ $0$ shares of risky asset (or stock) by time $T$.  He faces 
with the following risk/cost tradeoff:  if  he trades rapidly, this results in higher costs  for quickly executed orders and  market price impact; he can then  
split the order into several smaller blocks, but is then exposed to the risk of price depreciation during the trading horizon.  These liquidity effects  received  recently a considerable   interest  starting with the papers by Bertsimas and Lo \cite{berlo98}, and Almgren and Criss \cite{almcri01} in a discrete-time framework, and  further investigated among others in Obizhaeva and Wang \cite{obiwan05},  Schied an Sch\"oneborn 
\cite{schsch07},  or Rogers and Singh \cite{rogsin08} in a continuous-time model.  
These papers assume continuous trading with instantaneous trading rate inducing price impact.   In a continuous time 
market framework, we propose here a more realistic modeling by considering that trading takes place at discrete points in time  through an impulse control formulation,  and with a temporary price  impact depending  on the time interval between trades, and including a bid-ask spread.

We present the details of the model. Let $(\Omega,\Fc,\P)$ be a probability space equipped with a filtration $\F$ $=$ 
$(\Fc_{t})_{0\leq t\leq T}$ satisfying the usual conditions, and supporting a one dimensional Brownian motion $W$ on a finite horizon $[0,T]$, $T<\infty$. 
We denote by $P$ $=$ $(P_t)$ the market  price process  of the risky asset,  by  
$X_{t}$ the amount of money (or cash holdings),  by  $Y_{t}$ the number of shares in the stock held by the investor at time $t$, and by $\Theta_t$ the time interval between  time $t$ and the last trade before $t$.  We set $\R^*_+$ $=$ $(0,\infty)$ and $\R^*_-$ $=$ $(-\infty,0)$. 

 \vspace{2mm}

$\bullet$ \textit{Trading strategies.} We assume that the investor can only trade discretely on $[0,T]$. This is modelled through an impulse control strategy $\alpha$ $=$ $(\tau_{n},\zeta_{n})_{n\geq0}$: $\tau_{0}\leq\ldots\leq \tau_{n} \ldots\leq T$ are nondecreasing stopping times representing the trading times of the investor and 
$\zeta_{n}$, $n\geq 0$, are $\Fc_{\tau_{n}}-$measurable random variables valued in $\R$ and giving the number of stock purchased if 
$\zeta_{n}\geq0$ or selled if $\zeta_{n}<0$ at these times.  We denote by $\Ac$ the set of trading strategies. 
The sequence $(\tau_{n},\zeta_{n})$ may be a priori finite or infinite.  Notice also that we do not assume a priori that the sequence of 
trading times $(\tau_n)$ is strictly increasing. We  introduce the  lag variable tracking the time interval between trades: 
\beqs \label{Theta}
\Theta_t  &=& \inf\big\{ t - \tau_n:  \tau_n \leq t \}, \;\;\ t \in [0,T],
\enqs
which evolves according to
\beq \label{dynTheta}
\Theta_t  \; = \;   t - \tau_n, \;\;\; \tau_n \leq t < \tau_{n+1},  & & \Theta_{\tau_{n+1}} = 0, \;\;\; n \geq 0. 
\enq
The dynamics of the number of shares  invested in stock is given by:
\beq \label{Y}
Y_{t}  \;=  \;  Y_{\tau_{n}},\qquad \tau_{n}\leq t<\tau_{n+1},  & & 
Y_{\tau_{n+1}} \;  =  \;  Y_{\tau_{n+1}^-}+\zeta_{n+1}, \;\;\; n \geq 0. 
\enq 
 
 \vspace{2mm}

$\bullet$ \textit{Cost of illiquidity.}  The market price of the risky asset  process  follows a geometric Brownian motion:
\beq \label{dynP}
dP_t &=&  P_t (b dt + \sigma dW_t),
\enq
with constants $b$ and $\sigma$ $>$ $0$.  We do not consider a permanent price impact on the price, i.e. the lasting effect of large trader, but focus here on the effect of illiquidity, that is the price at which an investor will trade the asset.  Suppose now that the investor decides at time $t$ to make an order in stock shares of size  $e$. If the 
current market price is $p$, and the time lag from the last order is $\theta$, then the price he actually get  for the order $e$ is: 
\beq \label{Qimpact}
Q(e,p,\theta) & = &  p f (e,\theta),\enq
where  $f$ is a temporary price impact function from $\R\times [0,T]$ into $\R_{+}\cup\{\infty\}$. We assume that the Borelian function $f$ satisfies the following  liquidity and transaction cost properties: 

\vspace{1mm}

{\bf (H1f)} \hspace{7mm} $f(0,\theta)$ $=$ $1$, and $f(.,\theta)$ is nondecreasing for all $\theta\in[0,T]$, 
 
{\bf (H2f)} \hspace{6mm}   (i) $f(e,0)$ $=$ $0$ for $e$ $<$ $0$, and (ii) $f(e,0)$ $=$ $\infty$ for $e$ $>$ $0$, 
 
{\bf (H3f)} \hspace{6mm}  $\kappa_b$ $:=$ $\sup_{(e,\theta)\in\R^*_{-}\times[0,T]}f(e,\theta)$ $<$ $1$ and $\kappa_a$ $:=$ $\inf_{(e,\theta)\in\R_{+}^*\times[0,T]}f(e,\theta)$ $>$ $1$. 

\vspace{1mm}
 
\noindent Condition {\bf (H1f)}  means that  no trade incurs no impact on the market price, i.e. $Q(0,p,\theta)$ $=$ $p$,  
and a purchase (resp. a sale) of stock shares induces a cost (resp. gain) greater  (resp. smaller) than the market price, which increases (resp. decreases) with the size of the order. 
 In other words, we have $Q(e,p,\theta)$ $\geq$ (resp.  $\leq$) $p$ for $e$ $\geq$ (resp. $\leq$) $0$, and $Q(.,p,\theta)$ is nondecreasing.   
 Condition {\bf (H2f)}  expresses  the higher costs for immediacy in trading:  indeed, the immediate market resiliency is limited, and 
 the faster  the investor wants to liquidate  (resp. purchase) the asset, the deeper into the limit order book he will have to go, and lower 
 (resp. higher) will be the price for the shares of the asset sold  (resp. bought), with a zero (resp. infinite) limiting price for 
 immediate block sale (resp. purchase).  Condition {\bf (H2f)} also prevents the investor to pass orders at consecutive immediate times, which is 
  the case in  practice. Instead of  imposing a fixed arbitrary lag between orders, we shall see that condition {\bf (H2)} implies 
 that trading times are strictly increasing.   
 Condition {\bf (H3f)} captures a transaction cost effect:  at time $t$, $P_t$ is the market or mid-price, $\kappa_b P_t$ is the bid price, 
$\kappa_a P_t$ is the ask  price, and $(\kappa_a-\kappa_b)P_t$ is the bid-ask spead. We also assume some 
regularity  conditions on the temporary price impact function:

\vspace{1mm}

{\bf (Hcf)} \hspace{6mm}   (i) $f$ is continuous on $\R^*\times(0,T]$,  

\hspace{17mm} (ii)  $f$ is   $C^1$ on $\R_{-}^*\times[0,T]$ and $x$ $\mapsto$ $\Dth{f}$  
is bounded on $\R_{-}^*\times[0,T]$. 

\vspace{1mm}

\noindent 
A  usual  form (see e.g. \cite{lilfarman03}, \cite{potbou03}, \cite{almthuhau05}) of  temporary price impact and transaction cost function $f$, suggested by empirical studies is 
\beq \label{exf}
f(e,\theta) &=& e^{\lambda |\frac{e}{\theta}|^\beta {\rm sgn}(e)}  
\Big( \kappa_{a}\mathbf{1}_{e>0} +\mathbf{1}_{e=0}  +\kappa_{b}\mathbf{1}_{e<0} \Big),
\enq
with the convention $f(0,0)$ $=$ $1$.
Here $0<\kappa_b<1<\kappa_a$, $\kappa_a-\kappa_b$ is the bid-ask spread parameter,  $\lambda$ $>$ $0$  is the temporary price impact factor, and 
$\beta$ $>$ $0$  is the price impact exponent.   In our illiquidity modelling, we focus on the cost of  trading fast (that is the temporary price impact), and 
ignore as in Cetin, Jarrow and Protter \cite{cetjarpro04} and Rogers and Singh \cite{rogsin08}  the  permanent price impact of a large trade.  
This last effect could  be  included in our model, by assuming a jump of the price process at the trading date, depending on the order size,  
see e.g. He and Mamayski \cite{hemam04} and  Ly Vath, Mnif and Pham \cite{lyvmnipha07}.

\vspace{2mm}

$\bullet$ \textit{Cash holdings.} We assume a zero risk-free return, so that the bank account is constant between two trading times: 
\beq \label{Xconst}
X_t &=& X_{\tau_n}, \;\;\;\; \tau_n \leq t  <  \tau_{n+1}, \;\; n \geq 0.
\enq
When a discrete trading $\Delta Y_{t}$ $=$ $\zeta_{n+1}$ occurs at time $t$ $=$ $\tau_{n+1}$,  
this results  in a variation of the cash amount given by $\Delta X_t$ $:=$ $X_t-X_{t^-}$ $=$ $-\Delta Y_t . Q(\Delta Y_t,P_t,\Theta_{t^-})$ 
due to the illiquidity effects.  In other words,  we have
\beq
X_{\tau_{n+1}}  & = & X_{\tau_{n+1}^-} - \zeta_{n+1} Q(\zeta_{n+1},P_{\tau_{n+1}},\Theta_{\tau_{n+1}^-})  \nonumber  \\
&=&  X_{\tau_{n+1}^-} - \zeta_{n+1}  P_{\tau_{n+1}} f( \zeta_{n+1}, \tau_{n+1}-\tau_{n} ), \;\;\; n \geq 0.  \label{X}  
\enq
Notice that similarly as in the above cited papers dealing with continuous-time trading, 
we do not assume  fixed transaction fees to be paid at each trading.  
They are practically   insignificant with  respect to the price impact and bid-ask spread.   We can then not exclude a priori 
trading strategies with immediate trading times, i.e. $\Theta_{\tau_{n+1}^-}$ $=$ $\tau_{n+1} -\tau_{n}$ $=$ $0$ for some $n$.  
However, notice that  under condition {\bf (H2f)}, an immediate sale does not increase the cash holdings, i.e. 
$X_{\tau_{n+1}}$ $=$ $X_{\tau_{n+1}^-}$  $=$ $X_{\tau_n}$, 
while an immediate  purchase leads to a bankruptcy, i.e.  $X_{\tau_{n+1}}$ $=$ $-\infty$.

\vspace{2mm}

$\bullet$ \textit{Liquidation value and solvency constraint.}  A key issue in portfolio liquidation  is to define in an economically meaningful way what is the portfolio value of a position on cash and stocks. In our framework, we impose a no-short sale constraint on the trading strategies, i.e. 
\beqs \label{noshort}
Y_t &\geq & 0, \;\;\; 0 \leq t\leq T,
\enqs
which is in line with the  bank regulation following the financial crisis, 
and we consider  the liquidation function $L(x,y,p,\theta)$  representing the net wealth value 
that an investor with a cash amount $x$, would obtained by liquidating his stock position $y$ $\geq$ $0$ by a single block trade, when the market price is $p$ and given 
the time lag $\theta$ from  the last trade. It is  defined on $\R\times\R_+\times\R_+^*\times [0,T]$  by
\beqs
L(x,y,p,\theta) & = & x + ypf(-y,\theta), 
\enqs
and we impose the liquidation constraint on trading strategies: 
\beqs \label{liquidcon}
L(X_{t},Y_{t},P_{t},\Theta_{t}) & \geq & 0, \;\;\;  0 \leq t \leq T. 
\enqs   
We have $L(x,0,p,\theta)$ $=$ $x$, and under condition {\bf (H2f)}(ii), we notice that  $L(x,y,p,0)$ $=$ $x$ for $y$ $\geq$ $0$. 
We   naturally introduce the liquidation solvency region: 
\beqs
\Sc & = & \Big\{ (z,\theta)=(x,y,p,\theta) \in\R\times\R_+\times\R_+^*\times [0,T]:~y> 0~\mbox{ and }~L(z,\theta)> 0\Big\}. 
\enqs
We  denote its boundary and its closure by
\beqs
\partial \Sc & = & \partial_{y} \Sc\cup\partial_{L} \Sc~~\mbox{ and }~~\bar \Sc~=~\Sc \cup \partial \Sc,
\enqs
where 
\beqs
 \partial_{y} \Sc & = & \Big\{(z,\theta)=(x,y,p,\theta)\in\R\times\R_+\times\R_+^*\times [0,T]:~y= 0 \; \mbox{ and } \; x = L(z,\theta) \geq 0 \Big\}, \\
 \partial_{L} \Sc & = & \Big\{(z,\theta)=(x,y,p,\theta)  \in\R\times\R_+\times\R_+^*\times [0,T]:~ L(z,\theta)= 0\Big\}.
\enqs
We also denote by $D_0$ the corner line in $\partial\Sc$: 
\beqs
D_0 &=& \{0\}\times\{0\}\times\R_+^*\times [0,T] \; = \; \partial_y \Sc \cap \partial_L \Sc. 
\enqs

\begin{figure}
\begin{center}\label{Sc2d}
\includegraphics[angle=-90, width=16cm]
{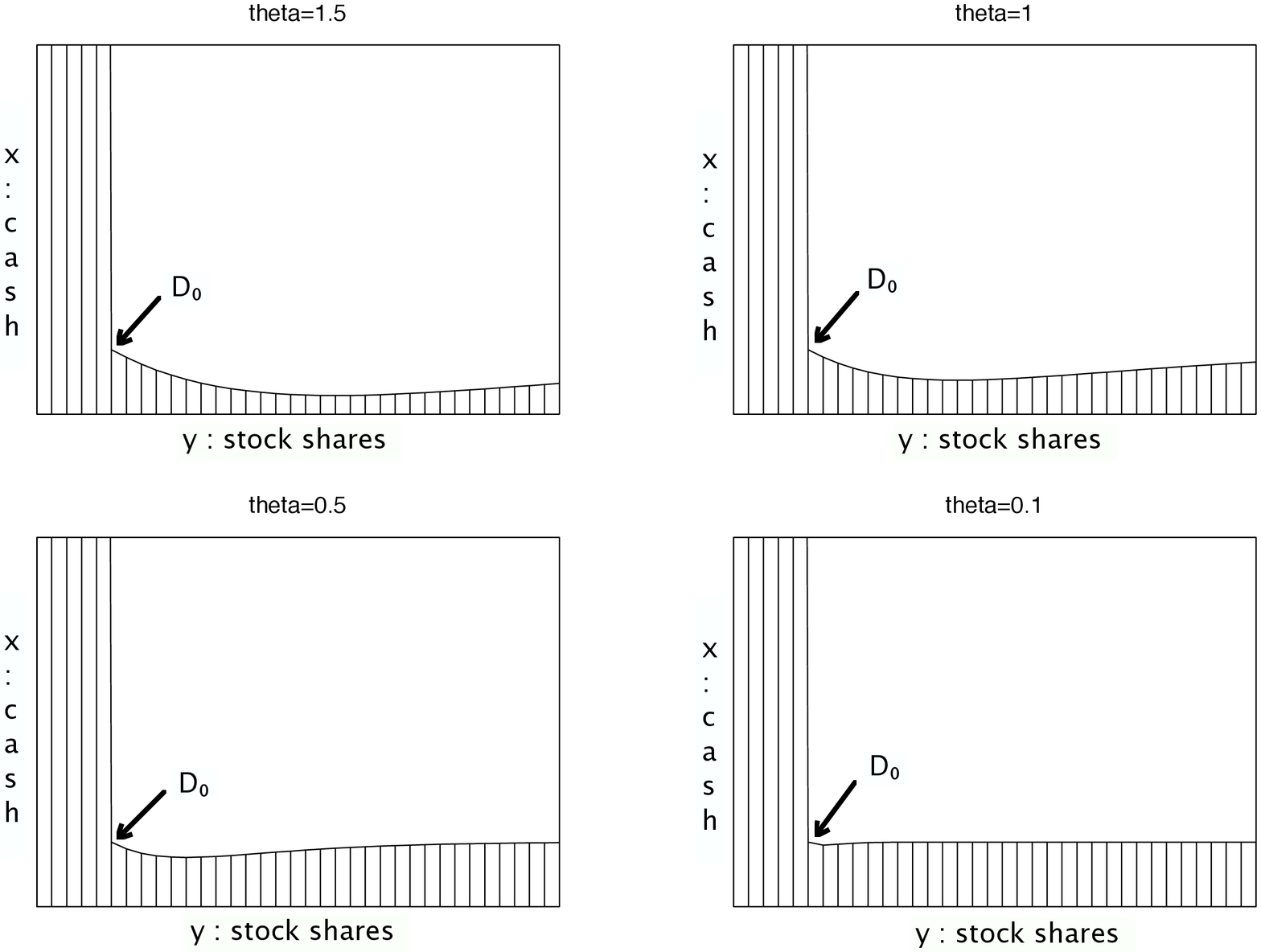}
\caption{Domain $\Sc$ in the nonhatched zone for fixed $p=1$ and $\theta$ evolving from $1.5$ to $0.1$. Here $\kappa_{b}=0.9$ and $f(e,\theta)=\kappa_{b}\exp(\frac{e}{\theta})$ for $e<0$. Notice that when $\theta$ goes to $0$, the domain converges to the open orthant 
$\R_{+}^*\times\R_{+}^*$.
}
\end{center}

\end{figure}

\begin{figure}
\label{S3d-theta-fix}
\begin{center}
\includegraphics[width=16cm,height=9cm]{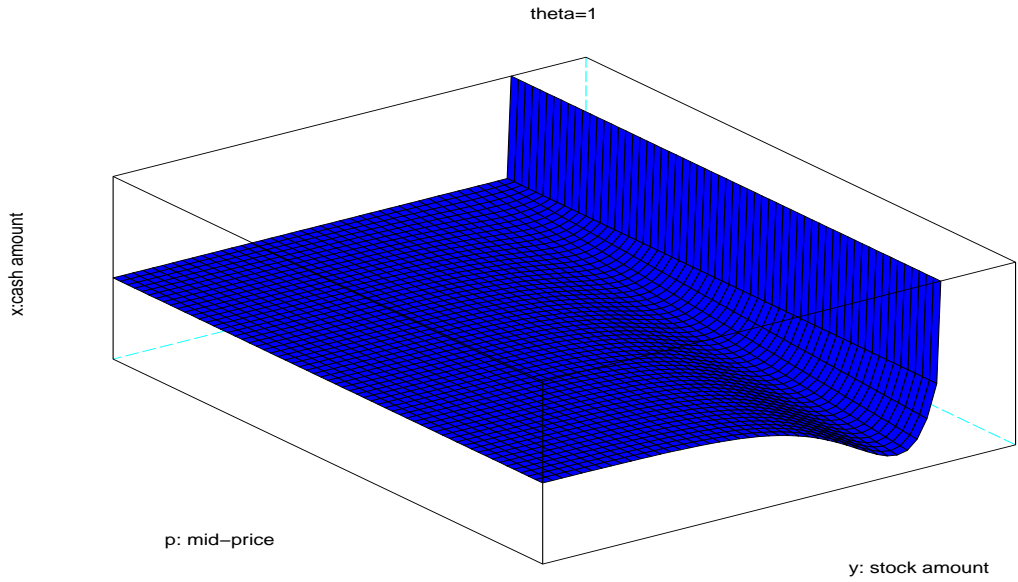}
\caption{Lower bound of the domain $\Sc$ for fixed  $\theta=1$. Here $\kappa_{b}=0.9$ and $f(e,\theta)=\kappa_{b}\exp(\frac{e}{\theta})$ for $e<0$. Notice that when $p$ is fixed, we obtain the Figure 1.
}
\end{center}
\end{figure}

\begin{figure}
\begin{center}\label{S3d-p-fix}
\includegraphics[width=15cm,height=9cm]{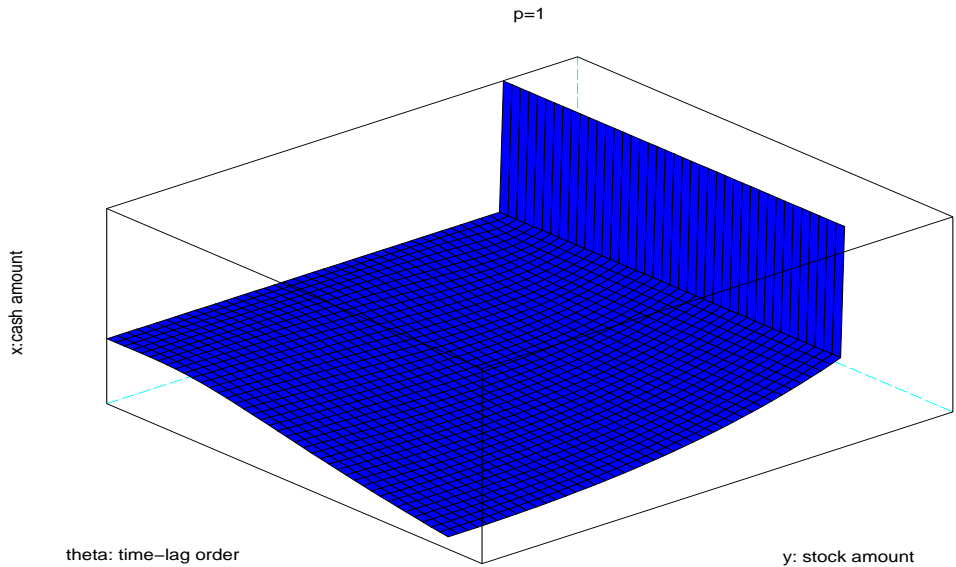}
\caption{Lower bound of the domain $\Sc$ for fixed $p=1$ with  $f(e,\theta)=\kappa_{b}\exp(\frac{e}{\theta})$ for $e<0$ and $\kappa_{b}=0.9$. Notice that when $\theta$ is fixed, we obtain the Figure 1. 
}
\end{center}
\end{figure}

\vspace{2mm}

$\bullet$   \textit{Admissible trading strategies.} Given  $(t,z,\theta)$ $\in$ $[0,T]\times\bar\Sc$,  we say that the impulse control strategy 
$\alpha$ $=$ $(\tau_{n},\zeta_{n})_{n\geq 0}$   is admissible, denoted by $\alpha$ $\in$ $\Ac(t,z,\theta)$,  
if $\tau_0$ $=$ $t-\theta$, $\tau_n$ $\geq$ $t$,  $n$ $\geq$ $1$,  and  the process 
$\{(Z_{s},\Theta_s) = (X_{s,}Y_{s},P_{s},\Theta_{s}), t\leq s\leq T\}$  solution to  \reff{dynTheta}-\reff{Y}-\reff{dynP}-\reff{Xconst}-\reff{X}, with  an initial state 
$(Z_{t^-},\Theta_{t^-})$ $=$ $(z,\theta)$ (and the convention that  $(Z_{t},\Theta_{t})$ $=$ $(z,\theta)$ if $\tau_1$ $>$ $t$),      
satisfies  $(Z_{s},\Theta_s)$ $\in$ $[0,T]\times\bar\Sc$ for all $s$ $\in$ $[t,T]$.  As usual, to  alleviate notations, we omitted 
 the dependence of $(Z,\Theta)$ in  $(t,z,\theta,\alpha)$, when there is no ambiguity.

\begin{Remark} \label{remadmi}
{\rm  Let  $(t,z,\theta)$ $\in$ $[0,T]\times\bar\Sc$, and consider the impulse control strategy $\alpha$ $=$ $(\tau_n,\zeta_n)_{n\geq 0}$, $\tau_0$ $=$ 
$t-\theta$,  consisting in liquidating immediately all the stock shares, and then doing no transaction anymore, i.e. $(\tau_1,\zeta_1)$ $=$ $(t,-y)$, and 
$\zeta_n$ $=$ $0$, $n$ $\geq$ $2$. The associated state process $(Z=(X,Y,P),\Theta)$ satisfies $X_s$ $=$ $L(z,\theta)$, $Y_s$ $=$ $0$,  
which shows that $L(Z_s,\Theta_s)$ $=$ $X_s$ $=$ $L(z,\theta)$ $\geq$ $0$, $t\leq s\leq T$,  and thus  $\alpha$ $\in$ $\Ac(t,z,\theta)$ $\neq$ $\emptyset$. 
}
\end{Remark}

\vspace{2mm}

$\bullet$ \textit{Portfolio liquidation problem.} We consider a utility function $U$ from $\R_+$ into $\R$, nondecreasing, concave, with  $U(0)$ $=$ $0$, and s.t.  there exists $K$ $\geq$ $0$ and $\gamma$ $\in$ $[0,1)$:
\beqs 
\hspace{-4.5cm} {\bf (HU)} \hspace{4cm} 0 \; \leq \; U(x) & \leq & K x^\gamma, \quad \forall x  \in\R_+.
\enqs 
The problem of optimal portfolio liquidation  is formulated as
\beq \label{defvliquid}
v(t,z,\theta) &=&   \sup_{\alpha\in\Ac_\ell(t,z,\theta)} \E\big[ U(X_T) \big], \;\;\; (t,z,\theta) \in [0,T]\times\bar\Sc,
\enq
where   $\Ac_\ell(t,z,\theta)$ $=$ $\big\{\alpha \in \Ac(t,z,\theta):~ Y_{T}~=~0 \big\}$ is nonempty by Remark \ref{remadmi}.  
Notice that for $\alpha$ $\in$ $\Ac_\ell(t,z,\theta)$, 
$X_T$ $=$ $L(Z_T,\Theta_T)$ $\geq$ $0$, so that the expectations in \reff{defvliquid}, and the value function $v$ 
are well-defined in $[0,\infty]$.  Moreover, by considering the particular strategy described in Remark \ref{remadmi},  which leads to 
a final liquidation value $X_T$ $=$ $L(z,\theta)$, we obtain a lower-bound for the value function;  
\beq \label{vlower}
v(t,z,\theta) & \geq & U(L(z,\theta)), \;\;\;  (t,z,\theta) \in [0,T]\times\bar\Sc.
\enq

\begin{Remark}
{\rm  We can shift  the terminal liquidation constraint in $\Ac_\ell(t,z,\theta)$ to a terminal liquidation utility  by considering 
the function $U_L$  defined on $\bar \Sc$ by: 
\beqs
U_{L}(z,\theta) & = & U(L(z,\theta)), \;\;\; (z,\theta) \in \bar\Sc. 
\enqs
Then, problem \reff{defvliquid} is written equivalently in 
\beq \label{defbarv}
\bar v(t,z,\theta) &=& \sup_{\alpha\in\Ac(t,z,\theta)} \E\Big[ U_L(Z_T,\Theta_T) \Big], \;\;\; (t,z,\theta) \in [0,T]\times\bar\Sc. 
\enq
Indeed, by observing that for all $\alpha$ $\in$ $\Ac_{\ell}(t,z,\theta)$,  we have 
$\E[U(X_T)]$ $=$ $\E[U_L(Z_T,\Theta_T)]$, and since $\Ac_\ell(t,z,\theta)$ $\subset$ $\Ac(t,z,\theta)$, it is clear that  $v$ $\leq$ $\bar v$.  Conversely, 
for any $\alpha$ $\in$ $\Ac(t,z,\theta)$ associated to the state controlled process $(Z,\Theta)$, 
consider the impulse control strategy $\tilde\alpha$ $=$ $\alpha$ $\cup$ $(T,-Y_T)$ consisting in 
liquidating all the stock shares $Y_T$  at time $T$.  The corresponding state process  $(\tilde Z,\tilde\Theta)$ satisfies clearly: 
 $(\tilde Z_s,\tilde\Theta_s)$ $=$ $(Z_s,\Theta_s)$ for $t\leq s<T$, and $\tilde X_T$ $=$ $L(Z_T,\Theta_T)$, $\tilde Y_T$ $=$ $0$, and so 
 $\tilde\alpha$ $\in$ $\Ac_\ell(t,z,\theta)$.  We deduce that $\E[U_L(Z_T,\Theta_T)]$ $=$ $\E[U(\tilde X_T)]$ $\leq$ $v(t,z,\theta)$, and so  by arbitrariness of  $\alpha$ in $\Ac(t,z,\theta)$,  $\bar v(t,z,\theta)$ $\leq$ $v(t,z,\theta)$. This proves the equality $v$ $=$ $\bar v$.   
 Actually, the above arguments also show that $\sup_{\alpha\in\Ac_\ell(t,z,\theta)} U(X_T)$ $=$ 
 $\sup_{\alpha\in\Ac(t,z,\theta)} U_L(Z_T,\Theta_T)$. 
}
\end{Remark}

\begin{Remark} \label{remcont}
{\rm  A continuous-time trading  version  of  our illiquid market model with stock price $P$ and temporary price impact $f$ 
 can be  formulated as follows.  The trading strategy  is  given by a 
$\F$-adated process $\eta$ $=$ $(\eta_t)_{0\leq t\leq T}$ representing the instantaneous trading rate, which means that the dynamics of the 
cumulated number of stock shares $Y$ is governed by: 
\beq \label{dynYcont}
dY_t &=& \eta_t dt. 
\enq
The cash holdings $X$ follows
\beq \label{dynXcont}
dX_t &=&  - \eta_t  P_t f(\eta_t) dt.
\enq
Notice that in a continuous-time trading formulation, the time interval between trades is $\Theta_t$ $=$ $0$ at any time $t$. Under condition 
{\bf (H2f)}, the liquidation value is then given at any time $t$ by: 
\beqs
L(X_t,Y_t,P_t,0) &=& X_t, \;\;\; 0 \leq t \leq T,
\enqs
and does not  capture  the position in stock shares, which  is economically not relevant. On the contrary, by explicitly considering the time interval between trades in our  discrete-time trading formulation, we take into account the position in stock. 
}
\end{Remark}

\section{Properties of the model}

\setcounter{equation}{0} \setcounter{Assumption}{0}
\setcounter{Theorem}{0} \setcounter{Proposition}{0}
\setcounter{Corollary}{0} \setcounter{Lemma}{0}
\setcounter{Definition}{0} \setcounter{Remark}{0}

In this section,  we show that the illiquid market model presented in the previous section displays some interesting and economically meaningful  properties on the admissible trading strategies and the optimal performance, i.e.  the value function.  
Let us consider  the impulse transaction function $\Gamma$  defined 
on $\R\times\R_+\times\R_+^*\times [0,T] \times\R$ into $\R\cup\{-\infty\}\times\R\times\R_+^*$ by:
\beqs
\Gamma(z,\theta,e) &=&  
\Big(x - e p f\big(e,\theta\big), y+ e,p\Big), 
\enqs
for $z$ $=$ $(x,y,p)$, and set $\bar\Gamma(z,\theta,e)$ $=$ $\big(\Gamma(z,\theta,e),0\big)$. 
This corresponds to the value of the state variable $(Z,\Theta)$ immediately after a  trading  at time $t$ $=$ $\tau_{n+1}$ 
of $\zeta_{n+1}$ shares of stock, i.e.  $(Z_{\tau_{n+1}},\Theta_{\tau_{n+1}})$ $=$ 
$\big(\Gamma(Z_{\tau_{n+1}^-},\Theta_{\tau_{n+1}^-},\zeta_{n+1}),0\big)$.   
We then define the set of admissible transactions: 
\beqs
\Cc(z,\theta) &=& \Big\{ e \in \R:  \big(\Gamma(z,\theta,e),0\big) \in \bar\Sc \Big\}, \;\;\;\;\; (z,\theta) \in \bar\Sc.
\enqs
This means that for any $\alpha$ $=$ $(\tau_n,\zeta_n)_{n\geq 0}$ $\in$ $\Ac(t,z,\theta)$ with associated state process $(Z,\Theta)$, we have 
$\zeta_{n}$ $\in$ $\Cc(Z_{\tau_{n}^-},\Theta_{\tau_n^-})$, $n$ $\geq$ $1$.
We define the impulse operator $\Hc$ by
\beqs
\Hc \varphi(t,z,\theta) &=& \sup_{e\in\Cc(z,\theta)} \varphi(t,\Gamma(z,\theta,e),0),  \;\;\;\;\; (t,z,\theta) \in [0,T]\times\bar\Sc.
\enqs

We also introduce the liquidation function of the (perfectly liquid) Merton model: 
\beqs
L_M(z) &=& x + py, \;\;\;\;\;  \forall z = (x,y,p) \in \R \times \R  \times \R_+^*. 
\enqs
For $(t,z=(x,y,p),\theta)$ $\in$ $[0,T]\times\bar\Sc$,  we denote by $(Z^{0,t,z},\Theta^{0,t,\theta})$  the state process starting from $(z,\theta)$ at time $t$, and without any  
impulse control strategy: it is given by  
\beqs
\Big(Z_s^{0,t,z},\Theta_s^{0,t,\theta}\Big) &=& (x,y,P_s^{t,p},\theta+ s-t), \;\;\; t \leq s \leq T, 
\enqs
where $P^{t,p}$ is the solution to \reff{dynP} starting from $p$ at time $t$.  Notice that  
$(Z^{0,t,z},\Theta^{0,t,\theta})$ is the continuous part of the state process $(Z,\Theta)$ controlled by  $\alpha$ $\in$ $\Ac(t,z,\theta)$. 
The infinitesimal generator $\Lc$ associated to the process 
$(Z^{0,t,z},\Theta^{0,t,\theta})$ is 
\beqs
\Lc \varphi + \Dth{\varphi}  &=& bp \Dp{\varphi} + \frac{1}{2} \sigma^2 p^2 \Dpp{\varphi} + \Dth{\varphi}. 
\enqs

\vspace{2mm}

We first prove a useful result on the set of admissible transactions.

\begin{Lemma} \label{lemsetadmi}
Assume that {\bf (H1f)}, {\bf (H2f)} and {\bf (H3f)}  hold. Then, 
for all $(z=(x,y,p),\theta)$ $\in$ $\bar\Sc$, the set $\Cc(z,\theta)$ is compact in $\R$ and satisfy 
\beq\label{inclC}
\Cc(z,\theta) &\subset& [-y,\bar e(z,\theta)],
\enq
where $-y\leq \bar e(z,\theta)<\infty$ is given by
\beqs
\bar e(z,\theta) &=& \left\{ 
				\begin{array}{cl}  
				\sup \Big\{ e \in\R:~ e p f(e,\theta) \leq x \Big\}\;, & \mbox{ if } \theta \; > \; 0 \\
				0\;, & \mbox{ if } \theta \; = \; 0.
				\end{array}
				\right.
\enqs
For $\theta$ $=$ $0$, \reff{inclC} becomes an equality : $\Cc(z,0)$ $=$ $[-y,0]$.

The set function  $\Cc$ is continous for the Hausdorff metric, i.e. if $(z_n,\theta_n)$ converges to $(z,\theta)$ in $\bar\Sc$, and 
$(e_n)$  is a sequence in $\Cc(z_n,\theta_n)$ converging to $e$, then $e$ $\in$ $\Cc(z,\theta)$.  Moreover, if $e\in\R$ $\mapsto$ $ef(e,\theta)$ is strictly increasing for $\theta\in(0,T]$, then for  $(z=(x,y,p),\theta)$ $\in$ $\partial_L\Sc$ with $\theta$ $>$ $0$, we have 
$\bar e(z,\theta)$ $=$ $-y$, i.e. $\Cc(z,\theta)$ $=$ $\{-y\}$. 
\end{Lemma}
{\bf Proof.} By definition of the impulse transaction function $\Gamma$ and the liquidation function $L$, we immediately see that the set of  admissible transactions is written as
\beq
\Cc(z,\theta) &=& \Big\{ e \in \R:~ x - e pf(e,\theta) \geq 0, \; \mbox{ and } \; y + e \geq 0 \Big\} \nonumber\\
&=& \Big\{ e\in\R:~ e p f(e,\theta) \leq x \Big\} \cap [-y,\infty) \; =: \; \Cc_1(z,\theta) \cap [-y,\infty). \label{egCC1} 
\enq
It is clear that $\Cc(z,\theta)$ is closed and bounded, thus a compact set. 
Under {\bf (H1f)} and {\bf (H3f)}, we have $\lim_{e\rightarrow\infty}epf(e,\theta)$ $=$ $\infty$. Hence we get $\bar e (z,\theta)$ $<$ $\infty$ and $\Cc_{1}(z,\theta)\subset(-\infty,\bar e (z,\theta)]$. From \reff{egCC1}, we get \reff{inclC}. Suppose $\theta$ $=$ $0$. Under {\bf (H2f)}, using $(z,\theta)$ $\in$ $\bar \Sc$, we have $\Cc_{1}(z,\theta)$ $=$ $\R_{-}$.  From \reff{egCC1}, we get $\Cc(z,\theta)$ $=$ $[-y,0]$.

Let us now  prove the continuity of the set of admissible transactions. Consider a sequence 
$(z_n=(x_n,y_n,p_n),\theta_n)$ in $\bar\Sc$ converging 
to $(z,\theta)$ $\in$ $\bar\Sc$, and  a sequence $(e_n)$  in  $\Cc(z_n,\theta_n)$  
converging to $e$. Suppose first that $\theta$ $>$ $0$. 
Then, for $n$ large enough, $\theta_n$ $>$ $0$ and  by observing that $(z,\theta,e)$ $\mapsto$ 
$\bar \Gamma (z,\theta,e)$ is continuous on  $\R\times\R_{+}\times\R_+^*\times\R_+^*\times\R$, we immediately deduce that 
$e$ $\in$ $\Cc(z,\theta)$. 
 In the case  $\theta$ $=$ $0$, writing $x_n-e_nf(e_n,\theta_n)$ $\geq$ $0$, using {\bf (H2f)}(ii) and sending $n$ to infinity, we see that $e$ should necessarily be nonpositive. By writing also that $y_n+e_n$ $\geq$ $0$, we get by sending $n$ to infinity that 
 $y+e$ $\geq$ $0$, and therefore $e$ $\in$ $\Cc(z,0)$ $=$ $[-y,0]$.

Suppose finally that $e\in\R$ $\mapsto$ $ef(e,\theta)$ is increasing, and fix  $(z=(x,y,p),\theta)$ $\in$ $\partial_L\Sc$, with $\theta$ $>$ $0$.  Then, 
$L(z,\theta)$ $=$ $0$, i.e. $x$ $=$ $-ypf(-y,\theta)$.  
Set $\bar e$ $=$ $\bar e(z,\theta)$.  By writing that $\bar e pf(\bar e,\theta)$ $\leq$ $x$ $=$ $-ypf(-y,\theta)$, and $\bar e$ $\geq$ $-y$, 
we deduce from the increasing monotonicity of $e$ $\mapsto$ $epf(e,\theta)$ that $\bar e$ $=$ $-y$. 
\ep

\begin{Remark} \label{remzetaneg}
{\rm The previous Lemma  implies in particular that $\Cc(z,0)$ $\subset$ $\R_-$, which means that an admissible transaction after an immediate trading should be necessarily a sale. In other words, given  $\alpha$ $=$ $(\tau_n,\zeta_n)_{n\geq 0}$ $\in$ $\Ac(t,z,\theta)$, $(t,z,\theta)$ $\in$ 
$[0,T]\times\bar\Sc$, if $\Theta_{\tau_n^-}$ $=$ $0$,  then $\zeta_n$ $\leq$ $0$.  
The continuity property of $\Cc$ ensures that the operator  $\Hc$ preserves the lower and upper-semicontinuity (see Appendix). 
This Lemma also asserts that, under the assumption of increasing monotonicity of $e$ $\rightarrow$ $ef(e,\theta)$, when the state is in the boundary 
$L$ $=$ $0$, then the only admissible transaction is to liquidate all stock shares.  This increasing monotonicity means that  the amount traded 
is increasing with the size of the order.  Such an assumption is satisfied in the example \reff{exf} of temporary price impact function $f$ 
for $\beta$ $=$ $2$, but is not fulfilled for $\beta$ $=$ $1$. In this case, the presence of illiquidity cost implies that it may be 
more advantageous to split the order size. 
}
\end{Remark}

\vspace{2mm}

We next  state  some useful  bounds  on the liquidation value associated to an admissible transaction.

\begin{Lemma} \label{lemliquidmerton}
Assume that  {\bf (H1f)} holds. Then, we have for all $(t,z,\theta)$ $\in$ $[0,T]\times\bar\Sc$:
\beq
0 \; \leq \; L(z,\theta) & \leq & L_M(z),  \label{LLM} \\
L_M(\Gamma(z,\theta,e)) & \leq & L_M(z), \;\;\; \forall e \in \R, \label{LMgam} \\
\sup_{\alpha\in\Ac(t,z,\theta)} L(Z_s,\Theta_s) & \leq & L_M(Z_s^{0,t,z}), \;\;\; t \leq s \leq T.  \label{LZ1}
\enq
Furthermore, under {\bf (H3f)}, we have for all $(z=(x,y,p),\theta)$ $\in$ $\bar\Sc$, 
\beq \label{LMgamstrict}
L_M(\Gamma(z,\theta,e)) & \leq & L_M(z) -  \min(\kappa_a-1,1-\kappa_b) |e| p, \;\;\; \forall e \in \R. 
\enq
\end{Lemma}
{\bf Proof.}
Under {\bf (H1f)}, we have $f(e,\theta)$ $\leq$ $1$ for all $e$ $\leq$ $0$,  which shows clearly  \reff{LLM}.  From the definition of $L_M$ and 
$\Gamma$, we see that for all $e$ $\in$ $\R$,
\beq \label{LMgamexpli}
L_M(\Gamma(z,\theta,e)) - L_M(z) &=&  e p \Big( 1 - f(e,\theta) \Big),
\enq
which yields   the inequality \reff{LMgam}.  Fix some arbitrary $\alpha$ $=$ $(\tau_n,\zeta_n)_{n\geq 0}$ $\in$ $\Ac(t,z,\theta)$ associated to the 
controlled state process $(Z,\Theta)$. When a transaction occurs at time $s$ $=$ $\tau_n$, $n$ $\geq$ $1$, the jump of 
$L_M(Z)$ is  nonpositive by \reff{LMgam}:
\beqs
\Delta L_M(Z_s) &=& L_M(Z_{\tau_n}) - L_M(Z_{\tau_n^-}) \; = \; L_M(\Gamma(Z_{\tau_n^-},\Theta_{\tau_n^-},\zeta_n)) - L_M(Z_{\tau_n^-}) 
\; \leq \; 0. 
\enqs
We deduce that the process $L_M(Z)$ is smaller than its continuous part equal to $L_M(Z^{0,t,z})$, and we then get  \reff{LZ1} with \reff{LLM}.  
Finally, under the additional condition {\bf (H3f)},  we easily obtain inequality \reff{LMgamstrict} from  relation \reff{LMgamexpli}. 
\ep

\vspace{2mm}

We now check that our liquidation problem is well-posed by stating a natural  upper-bound on the optimal performance, namely that  the value function  in our illiquid market model is bounded by the usual Merton bound in a perfectly liquid market.

\begin{Proposition}\label{PropBound}
Assume that {\bf (H1f)}  and  {\bf (HU)}  hold. 
Then, for all $(t,z,\theta)$ $\in$ $[0,T]\times\bar\Sc$,  the family 
$\{ U_L(Z_T,\Theta_T),  \alpha \in \Ac(t,z,\theta)\}$ is uniformly integrable, and  we have
\beq
v(t,z,\theta) \;  \leq \;  v_{0}(t,z) & := &\E\Big[ U\Big(  L_M \big( Z_{T}^{0,t,z} \big) \Big) \Big], \quad (t,z,\theta) \in[0,T]\times\bar\Sc, \nonumber \\
& \leq &  K e^{\rho(T-t)}  L_M(z)^\gamma,  \label{v0bound}
\enq
where $\rho$ is a positive constant s.t.
\beq\label{condrho}
\rho & \geq & \frac{\gamma}{1-\gamma} \frac{b^2}{2\sigma^2}.
\enq
\end{Proposition}
\textbf{Proof.}
From \reff{LZ1} and the nondecreasing monotonicity of $U$, we have for all $(t,z,\theta)$ $\in$ $[0,T]\times\bar\Sc$: 
\beqs
\sup_{\alpha\in\Ac_\ell(t,z,\theta)} U(X_T) \; = \; \sup_{\alpha\in\Ac(t,z,\theta)} U_L(Z_T,\Theta_T) & \leq & U(L_M(Z_T^{0,t,z})),
\enqs
and all the assertions of the Proposition will follow once we prove the inequality \reff{v0bound}. For this, 
consider  the nonnegative function $\varphi$ defined on $[0,T]\times\bar\Sc$ by:
\beqs
\varphi(t,z,\theta) & = & e^{\rho(T-t)} L_M(z)^\gamma~=~ e^{\rho(T-t)}\big( x+py \big)^\gamma,
\enqs
and notice that $\varphi$ is smooth $C^2$ on $[0,T]\times(\bar \Sc\setminus D_0)$.    
We claim that for $\rho$ $>$ $0$  large enough, the function $\varphi$ satisfies:
\beqs
-\Dt{\varphi} - \Dth{\varphi} - \Lc\varphi & \geq & 0, \;\;\;\;\; \mbox{ on }  \;  [0,T]\times(\bar \Sc\setminus D_0). 
\enqs 
Indeed, a straightforward calculation shows that for all $(t,z,\theta)\in[0,T]\times(\bar \Sc\setminus D_0)$:
\beq
& & -\Dt{\varphi}(t,z,\theta)-\Dth{\varphi}(t,z,\theta) - \Lc\varphi(t,z,\theta) \nonumber \\
& = &  e^{\rho(T-t)} L_M(z)^{\gamma-2} 
\Big[\Big( \sqrt{\rho} L_M(z)+\frac{b\gamma}{2\sqrt{\rho}}yp \Big)^2 
+ \Big(\frac{\gamma(1-\gamma)\sigma^2}{2}-\frac{b^2\gamma^2}{4\rho}\Big) y^2p^2\Big] \label{derivLM}
\enq   
which is nonegative under condition \reff{condrho}. 

Fix some $(t,z,\theta)\in[0,T]\times\bar \Sc$. If $(z,\theta)$ $=$ $(0,0,p,\theta)$ $\in$ $D_0$,  then we clearly have $v_0(t,z,\theta)$ $=$ $U(0)$, and inequality 
\reff{v0bound} is trivial. Otherwise, if $(z,\theta)$ $\in$ $\bar\Sc\setminus D_0$, then  the process $(Z^{0,t,z},\Theta^{0,t,\theta})$ satisfy $L_{M}(Z^{0,t,z},\Theta^{0,t,\theta})$ $>$ $0$. 
Indeed, Denote by $(\bar Z ^{t,z},\bar \Theta ^{t,\theta})$ the process starting from $(z,\theta)$ at $t$ and associated to the strategy consisting in liquidating all stock shares at $t$. Then we have $(\bar Z ^{t,z}_{s},\bar \Theta ^{t,\theta}_{s})$ $\in$ $\bar\Sc\setminus D_0$ for all $s\in[t,T]$ and hence $L_{M}(\bar Z ^{t,z}_{s},\bar \Theta ^{t,\theta}_{s})$$>$ $0$ for all $s\in[t,T]$.  Using \reff{LZ1} we get $L_{M}( Z ^{0,t,z}_{s}, \Theta ^{0,t,\theta}_{s})$ $\geq$ $L_{M}( \bar Z ^{t,z}_{s}, \bar  \Theta ^{t,\theta}_{s})$ $>$ $0$.

We can then  apply It\^o's formula to $\varphi(s,Z_s^{0,t,z},\Theta_s^{0,t,\theta})$ between  $t$ and  
$T_R$ $=$ $\inf\{ s\geq t:~|Z^{0,t,z}_{s}|\geq R \}\wedge T$:
\beqs
\E[\varphi(T_{R},Z^{0,t,z}_{T_{R}},\Theta^{0,t,\theta}_{T_R})] 
& = & \varphi(t,z)+\E\Big[\int_{t}^{T_{R}}\Big(\Dt{\varphi}+\Dth{\varphi}+ \Lc\varphi\Big)(s,Z^{0,t,z}_{s},\Theta_s^{0,t,\theta})ds\Big] \\
&\leq& \varphi(t,z).
\enqs
(The stochastic integral term vanishes in expectation since the integrand is bounded before $T_R$). 
By sending $R$ to infinity, we get by Fatou's lemma and since  $\varphi(T,z,\theta)$ $=$ $L_M(z)^\gamma$: 
\beqs
\E\Big[L_M(Z^{0,t,z}_{T})^\gamma\Big] & \leq & \varphi(t,z,\theta). 
\enqs
We conclude with the growth condition  {\bf (HU)}. 
\ep

\vspace{2mm}

As  a direct consequence of the previous proposition, we obtain the continuity of the value function on the boundary $\partial_y\Sc$, i.e. when we start with no stock shares.

\begin{Corollary} \label{coroSy}
Assume that {\bf (H1f)}  and  {\bf (HU)}  hold.  Then, the value function $v$ is con\-tinuous on $[0,T]\times\partial_y\Sc$, and we have
\beqs
v(t,z,\theta) &=& U(x), \;\;\; \forall t \in [0,T], (z,\theta)=(x,0,p,\theta) \in \partial_y\Sc. 
\enqs
In particular, we have $v(t,z,\theta)$ $=$ $U(0)$ $=$ $0$, for all $(t,z,\theta)$ $\in$ $[0,T]\times D_0$. 
\end{Corollary}
{\bf Proof.}  From the lower-bound \reff{vlower} and the upper-bound in Proposition \ref{PropBound}, we have for all 
$(t,z,\theta)$ $\in$ $[0,T]\times\bar\Sc$, 
\beqs
U\Big( x+ypf\big(-y,\theta\big) \Big) \; \leq \; v(t,z,\theta) & \leq & \E\big[U(L_M(Z_T^{0,t,z}))\big] \; = \;  \E\big[ U(x + y P_T^{t,p}) \big]. 
\enqs
These two inequalities imply  the required result.  
\ep

\vspace{2mm}

The following result states the finiteness of the total number of shares and amount traded.

\begin{Proposition}\label{CVseries}
Assume that {\bf (H1f)} and {\bf (H3f)} hold. Then, 
for any $\alpha$ $=$ $(\tau_n,\zeta_n)_{n\geq 0}$ $\in$ $\Ac(t,z,\theta)$, $(t,z,\theta)$ $\in$ $[0,T]\times\bar\Sc$,  we have
\beqs
\sum_{n\geq 1}|\zeta_{n}| ~ < ~ \infty, \;\;  
\sum_{n\geq 1}|\zeta_{n}| P_{\tau_{n}} \; < \;  \infty, &\mbox{ and }& 
\sum_{n\geq 1}|\zeta_{n}| P_{\tau_{n}}f\Big(\zeta_{n},\Theta_{\tau_{n}^-} \Big) ~ < ~ \infty, \;\;\;\;\; a.s.
\enqs 
\end{Proposition}
\textbf{Proof.} Fix  $(t,z=(x,y,p),\theta)$ $\in$ $[0,T]\times\bar\Sc$, and  $\alpha$ $=$ $(\tau_n,\zeta_n)_{n\geq 0}$ $\in$ $\Ac(t,z,\theta)$. 
Observe first  that the continuous part of the process $L_M(Z)$ is $L_M(Z^{0,t,z})$, and we denote its jump at time $\tau_n$ by 
$\Delta L_M(Z_{\tau_n})$ $=$ $L_M(Z_{\tau_n})-L_{M}(Z_{\tau_n^-})$. 
From the estimates \reff{LLM} and  \reff{LMgamstrict} in Lemma \ref{lemliquidmerton}, we then have almost surely for all $n$ $\geq$ $1$, 
\beqs
0  \; \leq \; L_M(Z_{\tau_n}) &=& L_M(Z^{0,t,z}_{\tau_n}) + \sum_{k=1}^n \Delta L_M(Z_{\tau_k})  \\
& \leq &  L_M(Z^{0,t,z}_{\tau_n})  -  \bar \kappa  \sum_{k=1}^n |\zeta_k| P_{\tau_k},
\enqs
where we set $\bar\kappa$ $=$ $\min(\kappa_a-1,1-\kappa_b)$ $>$ $0$.  We deduce  that for all $n$ $\geq$ $1$, 
\beqs
\sum_{k=1}^n |\zeta_k| P_{\tau_k} & \leq & \frac{1}{\bar\kappa} \sup_{s\in [t,T]} L_M(Z_s^{0,t,z}) \; = \;  
\frac{1}{\bar\kappa} \big( x + y \sup_{s\in [t,T]} P_s^{t,p}\big) \; < \; \infty, \;\;\; a.s.  
\enqs
This shows the almost sure convergence of the series $\sum_n  |\zeta_n| P_{\tau_n}$.  Moreover, since the price 
process $P$ is continous and strictly positive, we also obtain the convergence of the series 
$\sum_n |\zeta_n|$. 
Recalling that $f(e,\theta)$ $\leq$ $1$ for all $e$ $\leq$ $0$ and $\theta\in[0,T]$, we have for all $n$ $\geq$ $1$. 
\beq
\sum_{k=1}^n |\zeta_k| P_{\tau_k} f\big(\zeta_{k},\Theta_{\tau_{k}^-}\big) &=& 
\sum_{k=1}^n \zeta_k P_{\tau_k} f\big(\zeta_{k},\Theta_{\tau_{k}^-}\big)  
+  2 \sum_{k=1}^n |\zeta_k| P_{\tau_k} f\big(\zeta_{k},\Theta_{\tau_{k}^-}\big) \mathbf{1}_{\zeta_k  \leq 0} \nonumber \\
& \leq &  \sum_{k=1}^n \zeta_k P_{\tau_k} f\big(\zeta_{k},\Theta_{\tau_{k}^-}\big)  
+  2 \sum_{k=1}^n |\zeta_k| P_{\tau_k}.   \label{interserie}
\enq
On the other hand, we have 
\beqs
0~\leq~L_M(Z_{\tau_n})  & =  &  X_{\tau_n} + Y_{\tau_n} P_{\tau_n} \\
 & = &   x-\sum_{k=1}^n \zeta_{k}P_{\tau_{k}}f\big(\zeta_{k},\Theta_{\tau_{k}^-}\big)   
 +  (y +  \sum_{k=1}^n \zeta_{k})P_{\tau_{n}}.
\enqs
Together with \reff{interserie}, this implies that for all $n$ $\geq$ $1$, 
\beqs
\sum_{k=1}^n |\zeta_k| P_{\tau_k} f\big(\zeta_{k},\Theta_{\tau_{k}^-}\big) & \leq &  
x +  (y +  \sum_{k=1}^n|\zeta_{k}|) \sup_{s \in [t,T]} P_s^{t,p} +  2 \sum_{k=1}^n |\zeta_k| P_{\tau_k}. 
\enqs
The convergence of the series $\sum_n |\zeta_n| P_{\tau_n} f\big(\zeta_{n},\Theta_{\tau_{n}^-}\big)$ follows therefore from the convergence of 
the series  $\sum_n |\zeta_n|$ and $\sum_n  |\zeta_n| P_{\tau_n}$.
\ep

\vspace{2mm}

As a  consequence of the above results, we can now prove that in the optimal portfolio liquidation, 
it suffices to restrict to a finite number of  trading times,  which are  strictly increasing.  Given a  trading strategy 
$\alpha$ $=$ $(\tau_n,\zeta_n)_{n\geq 0}$ $\in$ $\Ac$, let us denote by $N(\alpha)$ 
the process counting the number of intervention times: 
\beqs
N_t(\alpha) &=& \sum_{n\geq 1}  \mathbf{1}_{\tau_{n} \leq  t}, \;\;\; 0 \leq t \leq T.
\enqs
We denote by $\Ac_\ell^b(t,z,\theta)$ the set of admissible trading strategies in $\Ac_\ell(t,z,\theta)$ with a 
finite number of  trading times, such that these trading times are  strictly increasing, namely: 
\beqs
\Ac_\ell^b(t,z,\theta) &=& \Big\{ \alpha = (\tau_n,\zeta_n)_{n\geq 0} \in \Ac_\ell(t,z,\theta):~N_T(\alpha) < \infty, \;\; \; a.s.  \\
& & \hspace{2.5cm} \; \mbox{ and } \;   \tau_n < \tau_{n+1} \;\; a.s., \; \;\;  0 \leq n \leq N_T(\alpha)-1 \Big\}. 
\enqs
For any $\alpha$ $=$ $(\tau_n,\zeta_n)_n$ $\in$ $\Ac_\ell^b(t,z,\theta)$, the associated state process $(Z,\Theta)$ satisfies 
$\Theta_{\tau_{n+1}^-}$ $>$ $0$, i.e. $(Z_{\tau_{n+1}^-},\Theta_{\tau_{n+1}^-})$ $\in$ $\bar\Sc^*$ $:=$  
$\Big\{ (z,\theta) \in \bar\Sc: \theta > 0 \Big\}$.  We also set $\partial_L\Sc^*$ $=$ $\partial_L\Sc$ $\cap$ $\bar\Sc^*$.

\begin{Theorem}  \label{prodis}
Assume that {\bf (H1f)}, {\bf (H2f)}, {\bf (H3f)}, {\bf (Hcf)}   and {\bf (HU)} hold. Then,   we have 
\beq \label{vfini}
v(t,z,\theta) &=& \sup_{\alpha\in\Ac_\ell^b(t,z,\theta)} \E\big[U(X_{T})\big], \;\;\; (t,z,\theta) \in [0,T]\times\bar\Sc.  
\enq 
Moreover,  we have
\beq \label{vfini2}
v(t,z,\theta) &=& \sup_{\alpha\in\Ac^b_{\ell_+}(t,z,\theta)} \E\big[U(X_{T})\big], \;\;\; (t,z,\theta) \in [0,T]\times(\bar\Sc\setminus\partial_L\Sc),  
\enq 
where $\Ac^b_{\ell_+}(t,z,\theta)$ $=$ $\{ \alpha \in \Ac_\ell^b(t,z,\theta): (Z_s,\Theta_s) \in (\bar\Sc\setminus\partial_L\Sc), t \leq s <  T\}$.
\end{Theorem}
\textbf{Proof.}
{\bf 1.}   Fix $(t,z,\theta)$ $\in$ $[0,T]\times\bar\Sc$, and  denote by $\bar\Ac_\ell^b(t,z,\theta)$ the set of admissible trading strategies in 
$\Ac_\ell(t,z,\theta)$ with a finite number of  trading times:
\beqs
\bar\Ac_\ell^b(t,z,\theta) &=& \Big\{ \alpha = (\tau_k,\zeta_k)_{k\geq 0} \in \Ac_\ell(t,z,\theta):~N_T(\alpha) \mbox{ is bounded a.s. } \Big\}.
\enqs
Given an arbitrary $\alpha$ $=$ $(\tau_k,\zeta_k)_{k\geq 0}$ $\in$ $\Ac_\ell(t,z,\theta)$ associated to the state process 
$(Z,\Theta)$ $=$ $(X,Y,P,\Theta)$,  let us consider the truncated trading strategy 
$\alpha^{(n)}$ $=$ $(\tau_k,\zeta_k)_{k\leq n}$ $\cup$ $(\tau_{n+1},-Y_{\tau_{n+1}^-})$,  which consists in liquidating all stock shares at time 
$\tau_{n+1}$. This strategy $\alpha^{(n)}$ lies in $\bar\Ac_\ell(t,z,\theta)$,  and is  
associated to the state process denoted by $(Z^{(n)},\Theta^{(n)})$.  We then have
\beqs
X_T^{(n)}  - X_T &=&   \sum_{k\geq n+1} \zeta_{k}P_{\tau_{k}}f\big(\zeta_{k},\Theta_{\tau_{k}^-}\big) \;  + \; 
Y_{\tau_{n+1}^-} P_{\tau_{n+1}} f\big(-Y_{\tau_{n+1}^-},\Theta_{\tau_{n+1}^-}\big). 
\enqs
Now,  from Proposition \ref{CVseries},  we have  
\beqs\label{cvseried2}
\sum_{k\geq n+1} \zeta_{k}P_{\tau_{k}}f\big(\zeta_{k},\Theta_{\tau_{k}^-}\big) & \longrightarrow & 0  \;\;\; a.s.  
\quad\mbox{ when  }\quad n\rightarrow \infty.
\enqs 
Moreover, since $0$ $\leq$ $Y_{\tau_{n+1}^-}$ $=$ $Y_{\tau_n}$ goes to $Y_T$ $=$ $0$ as $n$ goes to infinity,  by definition of $\alpha$ $\in$ 
$\Ac_\ell(t,z,\theta)$, and recalling that 
$f$ is smaller than $1$ on $\R_-\times[0,T]$, we deduce that 
\beqs
0 \; \leq \;  Y_{\tau_{n+1}^-} P_{\tau_{n+1}} f\big(-Y_{\tau_{n+1}^-}\Theta_{\tau_{n+1}^-}\big) & \leq & 
Y_{\tau_{n+1}^-} \sup_{s\in [t,T]} P_s^{t,p} \\
& \longrightarrow &   0  \;\;\; a.s.   \quad\mbox{ when  }\quad n\rightarrow \infty. 
\enqs
This proves that 
\beqs
X_T^{(n)} & \longrightarrow & X_T \;\;\; a.s.   \quad\mbox{ when  }\quad n\rightarrow \infty. 
\enqs
From Proposition \ref{PropBound}, the sequence $(U(X_T^{(n)}))_{n\geq 1}$  is uniformly integrable, and we can apply the dominated convergence theorem to get
\beqs
\E\big[U(X_{T}^{(n)})\big] & \longrightarrow & \E\big[U(X_{T})\big], \;\;\;  ~\mbox{ when  }~n\rightarrow\infty.  
\enqs
From the arbitrariness of $\alpha$ $\in$ $\Ac_\ell(t,z,\theta)$, this shows that 
\beqs
v(t,z,\theta) & \leq & \bar v^b(t,z,\theta) :=  \sup_{\alpha\in\bar\Ac_\ell^b(t,z,\theta)} \E\big[U(X_{T})\big],
\enqs
and actually the equality $v$ $=$ $\bar v^b$ since the other inequality $\bar v^b$ $\leq$ $v$  is trivial from the inclusion $\bar\Ac_\ell^b(t,z,\theta)$ $\subset$ 
$\Ac_\ell(t,z,\theta)$.

\noindent \textbf{2.} Denote by $v^b$  the value function in the r.h.s. of \reff{vfini}. It is clear that $v^b$ $\leq$  $\bar v^b$ $=$ $v$ since 
$\Ac_\ell^b(t,z,\theta)$ $\subset$ $\bar\Ac_\ell^b(t,z,\theta)$.  
To prove the reverse inequality we need first to study the behavior of optimal strategies at time $T$. Introduce the set 
\beqs
\tilde \Ac ^b_{\ell}(t,z,\theta) & = & \Big\{\alpha = (\tau_{k},\zeta_{k})_{k}\in\Ac ^b_{\ell}(t,z,\theta) ~:~\#\{k~:~\tau_{k}=T\}\leq 1\Big\}, 
\enqs
and denote by $\tilde v ^b$ the associated value function. 
Then we have $\tilde v ^b$ $\leq$ $\bar v ^b$.  Indeed, let $\alpha$ $=$ $(\tau_{k},\zeta_{k})_k$ be some arbitrary element in  
$\bar\Ac_{\ell}^b(t,z,\theta)$,  $(t,z=(x,y,p),\theta)$ $\in$ $[0,T]\times\bar\Sc$. If $\alpha\in\tilde \Ac_\ell^b(t,z,\theta)$ then we have $\tilde v ^{b} (t,z,\theta)\geq \E\Big[U_{L}(Z_{T},\Theta_{T})\Big]$, where $(Z,\Theta)$
  denotes the process associated to $\alpha$. 
  Suppose now that $\alpha\notin\tilde \Ac_\ell^b(t,z,\theta)$. Set $m$ $=$ $\max\{k~:~\tau_{k}<T\}$. Then define the stopping time $\tau':=\frac{\tau_{m}+T}{2}$ and the $\Fc_{\tau'}$-measurable random variable   
$\zeta'$ $:=$ arg$\max\{ ef(e,T-\tau_{m})~:~e\geq -Y_{\tau_{m}} \}$. Define the strategy $\alpha'$ $=$ $(\tau_{k},\zeta_{k})_{k\leq m}\cup(\tau', Y_{\tau_{m}}-\zeta')\cup(T,\zeta')$. From the construction of $\alpha'$, we easily check that $\alpha'\in \tilde \Ac^b(t,z,\theta)$ and $\E\Big[U_{L}(Z_{T},\Theta_{T})\Big]$ $\leq$ $\E\Big[U_{L}(Z'_{T},\Theta'_{T})\Big]$ where $(Z',\Theta')$
  denotes the process associated to $\alpha'$. Hence, we get $\tilde v ^b$ $\geq$ $\bar v ^b$.

We now prove that $v^b$ $\geq$ $\tilde v ^b$. Let  $\alpha$ $=$ $(\tau_{k},\zeta_{k})_k$ be some arbitrary element in  
$\tilde\Ac_{\ell}^b(t,z,\theta)$,  $(t,z=(x,y,p),\theta)$ $\in$ $[0,T]\times\bar\Sc$. 
 Denote by $N$ $=$ $N_T(\alpha)$ the a.s. finite number of trading times in $\alpha$.  
We set  $m$ $=$ $\inf\{  0 \leq k \leq N-1:~\tau_{k+1} = \tau_k \}$ and  $M$ $=$ $\sup\{ m+1 \leq k\leq N:~\tau_{k} = \tau_m\}$ with the convention 
that $\inf\emptyset$ $=$ $\sup\emptyset$ $=$ $N+1$. 
We then define 
$\alpha'$ $=$ $(\tau_k',\zeta_k')_{0\leq k\leq N-(M-m)+1}$ $\in$ $\Ac$ by:
\beqs
(\tau_k',\zeta_k') &=& \left\{
				\begin{array}{cl} 
				(\tau_k,\zeta_k), &   \mbox{ for } \; 0 \leq k < m   \\
				(\tau_m=\tau_M,\sum_{k=m}^M \zeta_k), & \mbox{ for } \;  k=m \mbox{ and } \; m < N, \\
				(\tau_{k+M-m},\zeta_{k+M-m}), &  \mbox{ for } \;  m+1 \leq k\leq N-(M-m) \mbox{ and } \; m < N,\\
				(\tau',\sum_{l=m+1}^M\zeta_{l}) & \mbox{ for } k=N-(M-m)+1\; 
				\end{array}
				\right. 
\enqs
where $\tau'=\frac{\hat \tau+T}{2}$ with $\hat \tau$ $=$ $\max\{\tau_{k}~:~\tau_{k}<T\}$,
and we denote by $(Z'=(X',Y',P),\Theta')$  the associated state process.   
It is clear that  $(Z_s',\Theta_s')$ $=$ $(Z_s,\Theta_s)$ for  $t\leq s<\tau_m$, and so 
$X'_{(\tau)'^-}$ $=$ $X_{(\tau')^-}$, $\Theta'_{(\tau')^-}$ $=$  $\Theta_{(\tau')^-}$.  Moreover, since 
$\tau_m$ $=$ $\tau_M$, we have $\Theta_{\tau_k^-}$ $=$ $0$ for $m+1\leq k\leq M$.  From Lemma \ref{lemsetadmi} 
(or Remark \ref{remzetaneg}),  this implies that $\zeta_k$ $\leq$ $0$ for $m+1\leq k\leq M$, and  so 
$\zeta_{N-(M-m)+1}'$  $=$ $\sum_{k=m+1}^M \zeta_k$ $\leq$ $0$.  We also recall that immediate sales  does not increase the cash holdings, so that 
$X_{\tau_k}$ $=$ $X_{\tau_m}$ for  $m+1\leq k\leq M$.  
We then get
\beqs
X'_{T} &=& X_{T} -  \zeta_{N-(M-m)+1}' P_{\tau'} f\big(\zeta_{N-(M-m)+1}',\Theta'_{(\tau')^-}\big) \\
& \geq & X_{T}. 
\enqs
 Moreover, we have $Y_T'$ $=$ $y+\sum_{k=1}^N \zeta_k$ $=$ $Y_T$ $=$ $0$.  
 By construction, notice that  $\tau_0'<\ldots<\tau_{m+1}'$.  Given an arbitrary   $\alpha$ $\in$ 
$\bar\Ac_\ell^b(t,z,\theta)$, we can then construct by induction a trading strategy $\alpha'$ $\in$ $\Ac_\ell^b(t,z,\theta)$ such that 
$X_T'$ $\geq$ $X_T$ a.s.  By the nondecreasing monotonicity of the utility function $U$, this yields
\beqs
\E[ U(X_T)] & \leq & \E[U(X_T')] \; \leq \; v^b(t,z,\theta),
\enqs
 and we conclude from the arbitrariness of $\alpha$ $\in$ $\tilde \Ac_\ell^b(t,z,\theta)$: $\tilde v^b$ $\leq$ $v^b$, and thus  
 $v$ $=$ $\bar v^b$ $=$ $\tilde v^b$ $=$ $v^b$.

\noindent \textbf{3.}   Fix now  an element $(t,z,\theta)$ $\in$ $[0,T]\times(\bar\Sc\setminus\partial_L\Sc)$, and denote by $v_+$ the r.h.s of 
\reff{vfini2}. It is clear that $v$ $\geq$ $v_+$. Conversely, take some arbitrary $\alpha$ $=$ $(\tau_k,\zeta_k)_k$ $\in$ 
$\Ac_\ell^b(t,z,\theta)$, associated with the state process $(Z,\Theta)$, and denote by $N$ $=$ $N_T(\alpha)$ the finite number of trading times in 
$\alpha$. Consider the first time  before $T$ when the liquidation value reaches zero, i.e. $\tau^\alpha$ $=$ 
$\inf\{t\leq s\leq T: L(Z_s,\Theta_s) = 0\}\wedge T$ with the convention $\inf\emptyset$ $=$ $\infty$. We claim that there exists 
$1\leq m\leq N+1$ (depending on $\omega$ and $\alpha$) such that $\tau^\alpha$ $=$ $\tau_m$,  with the convention that $m$ $=$ 
$N+1$, $\tau_{N+1}$ $=$ $T$  if $\tau^\alpha$ $=$ $T$. 
On the contrary,  there would  exist $1\leq k\leq N$ such that  $\tau_k < \tau^\alpha<\tau_{k+1}$,  
and $L(Z_{\tau^\alpha},\Theta_{\tau^\alpha})$ $=$ $0$. Between $\tau_k$ and $\tau_{k+1}$, there is no trading, and so 
$(X_s,Y_s)$ $=$ $(X_{\tau_k},Y_{\tau_k})$, $\Theta_s$ $=$ $s-\tau_k$ for $\tau_k\leq s<\tau_{k+1}$. We then get 
\beq \label{LZ}
L(Z_s,\Theta_s) &=& X_{\tau_k} + Y_{\tau_k} P_s f\big(-Y_{\tau_k},s-\tau_k\big), \;\; \tau_k \leq s < \tau_{k+1}. 
\enq
Moreover, since  $0$ $<$ $L(Z_{\tau_k},\Theta_{\tau_k})$ $=$ $X_{\tau_k}$,  and  $L(Z_{\tau^\alpha},\Theta_{\tau^\alpha})$ $=$ $0$, we see with 
\reff{LZ} for $s$ $=$ $\tau^\alpha$ that 
$Y_{\tau_k}P_{\tau^{\alpha}}f\big(-Y_{\tau_k},\tau^{\alpha}-\tau_{k}\big)$ should necessarily be strictly negative: $Y_{\tau_k}P_{\tau^{\alpha}}f\big(-Y_{\tau_k},\tau^{\alpha}-\tau_{k}\big)$ $<$ $0$, a contradiction with the admissibility conditions and the nonnegative property of $f$.

We then have $\tau^\alpha$ $=$ $\tau_m$ for some $1\leq m\leq N+1$.  Observe that if  $m$ $\leq$ $N$,  i.e. $L(Z_{\tau_m},\Theta_{\tau_m})$ $=$ $0$,  then $U(L(Z_T,\Theta_T))$ $=$ $0$. 
Indeed, suppose that $Y_{\tau_m}$ $>$ $0$ and $m$ $\leq$ $N$. 
From the admissibility condition, and by  It\^o's formula to $L(Z,\Theta)$ in \reff{LZ} between $\tau^\alpha$ and $\tau_{m+1}^-$, we get
\beq
0 \; \leq \; L(Z_{\tau_{m+1}^-},\Theta_{\tau_{k+1}^-}) &=&  L(Z_{\tau_{m+1}^-},\Theta_{\tau_{m+1}^-}) - L(Z_{\tau^\alpha},\Theta_{\tau^\alpha}) 
\nonumber\\
&=&  \int_{\tau^\alpha}^{\tau_{m+1}}  Y_{\tau_m} P_s \Big[ \beta(Y_{\tau_m},s-\tau_m) ds + \sigma  f\big(-Y_{\tau_k},s-\tau_m\big)  dW_s \Big],~~\qquad
 \label{Lalpha}
\enq
where $\beta(y,\theta)$ $=$ $b f(-y,\theta) + \Dth{f}(-y,\theta)$ is bounded on $\R_+\times[0,T]$ by {\bf (Hcf)}(ii).  
Since the integrand in the above stochastic integral w.r.t Brownian motion $W$  is strictly positive, thus nonzero,  we must have $\tau^\alpha$ $=$ 
$\tau_{m+1}$. Otherwise, there is a nonzero probability that the r.h.s. of \reff{Lalpha} becomes strictly negative, a contradiction with the inequality 
\reff{Lalpha}.

Hence we get 
$Y_{\tau_m}$ $=$ $0$, and thus $L(Z_{\tau_{m+1}^-},\Theta_{\tau_{m+1}^-})$ $=$ $X_{\tau_m}$ $=$ $0$. From the Markov feature of the  model and Corollary \ref{coroSy}, we then have
\beqs
\E\Big[U\Big(L(Z_{T},\Theta_{T})\Big)\Big|\Fc_{\tau_{m}}\Big] & \leq & v(\tau_{m},Z_{\tau_{m}},\Theta_{\tau_{m}})~=~U(X_{\tau_{m}}) ~=~0. 
\enqs
Since $U$ is nonnegative, this implies that $U\big(L(Z_{T},\Theta_{T})\big)$ $=$ $0$. 
 Let us next consider the trading strategy $\alpha'$ $=$ $(\tau_k',\zeta_k')_{0\leq k\leq (m-1)}$ 
 $\in$ $\Ac$ 
consisting in following $\alpha$ until time $\tau^\alpha$, and 
liquidating all stock shares at time $\tau^\alpha$ $=$ $\tau_{m-1}$,
and defined by:
\beqs
(\tau_k',\zeta_k') &=& \left\{
				\begin{array}{cl} 
				(\tau_k,\zeta_k), &   \mbox{ for } \; 0 \leq k < m-1   \\
				\Big(\tau_{m-1},-Y_{\tau_{(m-1)}^-}\Big), & \mbox{ for } \;  k=m-1, 
				\end{array}
				\right. 
\enqs
and we denote by $(Z',\Theta')$ the associated state process. It is clear that $(Z'_s,\Theta_s')$ $=$ $(Z_s,\Theta_s)$ for 
$t\leq s<\tau_{m-1}$, and so $L(Z_s',\Theta_s')$ $=$ $L(Z_s,\Theta_s)$  $>$ $0$ for $t\leq s\leq\tau_{m-1}$. 
The liquidation at time $\tau_{m-1}$ (for $m$ $\leq$ $N$) yields 
$X_{\tau_{m-1}}$ $=$ $L(Z_{\tau_{m-1}^-},\Theta_{\tau_{m-1}^-})$ $>$ $0$, and $Y_{\tau_{m-1}}$ $=$ $0$.  Since there is no more trading after time $\tau_{m-1}$, 
the liquidation value for $\tau_{m-1}\leq s \leq T$ is given by: $L(Z_s,\Theta_s)$ $=$ $X_{\tau_{m-1}}$ $>$ $0$.  This shows that $\alpha'$ $\in$ 
$\Ac^b_{\ell_+}(t,z,\theta)$.  When $m$ $=$ $N+1$, we have $\alpha$ $=$ $\alpha'$, and so $X_T'$ $=$ $L(Z_T',\Theta_T')$ $=$ 
$L(Z_T,\Theta_T)$ $=$ $X_T$. 
For $m$ $\leq$ $N$, we have  $U(X_T')$ $=$ $U(L(Z_T',\Theta_T'))$ $\geq$ $0$ $=$ $U(L(Z_T,\Theta_T))$ $=$ $U(X_T)$. 
We then get $U(X_T')$ $\geq$ $U(X_T)$ a.s., and so
\beqs
\E[ U(X_T)] & \leq & \E[U(X_T')] \; \leq \; v_+(t,z,\theta).
\enqs
We conclude from the arbitrariness of $\alpha$ $\in$ $\bar\Ac_\ell^b(t,z,\theta)$: $v$ $\leq$ $v_+$, and thus  
$v$ $=$ $v_+$.
\ep

\vspace{2mm}
\begin{Remark} 
{\rm If we suppose that the function $e\in\R$ $\mapsto$ $ef(e,\theta)$ is increasing for $\theta\in(0,T]$, we get the value of $v$ on the bound $\partial_L\Sc^*$: 
$v(t,z,\theta)$ $=$ $U(0)$ $=$ $0$ for $(t,z=(x,y,p),\theta) \in [0,T]\times\partial_L\Sc^*$. Indeed, fix some point $(t,z=(x,y,p),\theta) \in [0,T]\times\partial_L\Sc^*$, and consider an arbitrary $\alpha$ $=$ $(\tau_k,\zeta_k)_k$ $\in$ $\Ac_{\ell}^b(t,z,\theta)$ with state process $(Z,\Theta)$, and denote by $N$ the number of trading times. 
We distinguish two cases:  (i) If $\tau_1$ $=$ $t$, then by Lemma \ref{lemsetadmi}, the transaction $\zeta_1$ is equal to $-y$, which leads to 
$Y_{\tau_1}$ $=$ $0$, and 
 a liquidation value $L(Z_{\tau_1},\Theta_{\tau_1})$ $=$ $X_{\tau_1}$ $=$ $L(z,\theta)$ $=$ $0$.  At the next trading date $\tau_2$ (if it exists), 
 we get $X_{\tau_2^-}$ $=$ $Y_{\tau_2^-}$ $=$ $0$ with liquidation value $L(Z_{\tau_2^-},\Theta_{\tau_2^-})$ $=$ $0$, and by using again 
 Lemma \ref{lemsetadmi}, we see that after the transaction at $\tau_2$, we shall also obtain $X_{\tau_2}$ $=$ $Y_{\tau_2}$ $=$ $0$. By induction, 
 this leads at the final trading time to $X_{\tau_N}$ $=$ $Y_{\tau_N}$ $=$ $0$, and finally to $X_T$ $=$ $Y_T$ $=$ $0$. 
(ii) If $\tau_1$ $>$ $t$, we claim that  $y$ $=$ $0$.  On the contrary,  by arguing similarly as in \reff{Lalpha} between $t$ and $\tau_1^-$, 
we have then proved that any admissible trading strategy $\alpha$  $\in$ $\Ac_{\ell}^b(t,z,\theta)$ provides  a final liquidation value $X_T$ $=$ $0$, 
and so 
\beq \label{vL=0}
 v(t,z,\theta) &=& U(0) = 0, \;\;\;\;\;  \forall (t,z,\theta) \in [0,T]\times\partial_L\Sc^*. 
 \enq
}
\end{Remark}

\vspace{2mm}

\begin{Remark} \label{remdiscont}
{\rm The representation \reff{vfini} of the optimal portfolio liquidation reveals inte\-resting  economical and mathematical  features.  
It shows  that the liquidation problem in a continuous-time  illiquid market model  with discrete-time orders and temporary price impact 
with the presence of a bid-ask spread as considered in this paper, 
leads to nearly optimal trading strategies with a finite number  of  orders  and with strictly increasing trading times. 
While most models dealing with  trading strategies  via an impulse control formulation assumed fixed transaction fees in order to justify  
the discrete nature of trading times, we prove in this paper that discrete-time  trading  appears  naturally as a  consequence 
of temporary price impact and bid-ask spread. 

The representation \reff{vfini2} shows that when we are in an initial state with strictly positive liquidation value, then we can restrict in the 
optimal portfolio liquidation problem to admissible trading strategies with strictly positive liquidation value up to time $T^-$. The relation 
 \reff{vL=0} means that when the initial state has a zero liquidation value, which is not a result of an  immediate trading time, then the liquidation 
 value will stay at zero until the final horizon.   
}
\end{Remark}

\section{Dynamic programming and viscosity properties}

\setcounter{equation}{0} \setcounter{Assumption}{0}
\setcounter{Theorem}{0} \setcounter{Proposition}{0}
\setcounter{Corollary}{0} \setcounter{Lemma}{0}
\setcounter{Definition}{0} \setcounter{Remark}{0}

In the sequel,  the conditions  {\bf (H1f)}, {\bf (H2f)}, {\bf (H3f)}, {\bf (Hcf)}   and {\bf (HU)} stand in force, and are not recalled in the statement of Theorems and Propositions.

We use a dynamic programming approach to derive the equation satisfied by the value function of our optimal portfolio liquidation problem. 
Dynamic programming principle (DPP)  for impulse controls  was frequently used  starting from the works by Bensoussan and Lions \cite{benlio82}, 
and then considered e.g. in \cite{tanyon93},  \cite{okssul06}, \cite{lyvmnipha07} or \cite{sey08}. In our context 
(recall the expression \reff{defbarv} of the value function), this is formulated as: 

\vspace{2mm}

\noindent \textsc{Dynamic programming principle (DPP).} For all $(t,z,\theta)$ $\in$ $[0,T]\times\bar\Sc$, we have 
\beq\label{DPPgen}
v(t,z,\theta) & = & \sup_{\alpha\in\Ac(t,z,\theta)} \E[v(\tau,Z_{\tau},\Theta_\tau)],
\enq
where $\tau$ $=$ $\tau(\alpha)$ is any stopping time valued in $[t,T]$ eventually depending on the strategy $\alpha$ in \reff{DPPgen}.
More precisely we have :  
\begin{enumerate}[(i)]
\item for all $\alpha$ $\in$ $\Ac(t,z,\theta)$, for all $\tau$ $\in$ $\Tc_{t,T}$, the set of stopping times valued in $[t,T]$:
\beq\label{DPP1}
 \E[v(\tau,Z_{\tau},\Theta_\tau)] & \leq &  v(t,z,\theta)
\enq 

\item for all $\varepsilon$ $>$ $0$, there exists $\hat \alpha ^\varepsilon$ $\in$ $\Ac(t,z,\theta)$ s.t. for all $\tau$ $\in$ $\Tc_{t,T}$: 
\beq
\label{DPP2}
 v(t,z,\theta)  -  \varepsilon  & \leq & \E[v(\tau,\hat Z^\varepsilon_{\tau},\hat \Theta_\tau^\eps)],
\enq
with $(\hat Z^\varepsilon,\hat\Theta^\eps)$ the state process controlled by $\hat \alpha^\varepsilon$.
\end{enumerate}

\vspace{2mm}

The corresponding dynamic programming Hamilton-Jacobi-Bellman (HJB)  
equation is a quasi-variational inequality (QVI) written as: 
\beq\label{QVIv}
\min\Big[ -\Dt{v} -\Dth{v} - \Lc v~,~v -\Hc v \Big] & = & 0, \;\;\;\; \mbox{ in } [0,T)\times\bar\Sc,
\enq
together with the relaxed terminal condition: 
\beq \label{condterm}
\min\big[ v  - U_L~,~v - \Hc v \big] &=& 0,  \;\;\;\; \mbox{ in }  \{T\} \times\bar\Sc.
\enq

The rigorous derivation of  the HJB equation satisfied by the value function from the dynamic programming principle is achieved by means of the notion of viscosity solutions, and is by now rather classical in the modern approach of stochastic control (see e.g. 
the books  \cite{fleson93} and \cite{pha09}).  There are some specificities  here related to the impulse control and the liquidation state constraint, 
and we recall in Appendix,  definitions of (discontinuous) constrained viscosity solutions for parabolic QVIs.  
The main result of this section is stated as follows.

\begin{Theorem}
The value function $v$ is a constrained viscosity solution to \reff{QVIv}-\reff{condterm}.
\end{Theorem} 
{\bf Proof.}
The proof of the viscosity supersolution property on $[0,T)\times\Sc$ and the viscosity subsolution property on $[0,T)\times\bar\Sc$ 
follows the same lines of arguments as in \cite{lyvmnipha07}, and is then omitted here. We focus on the terminal condition \reff{condterm}.  

We first check the viscosity supersolution property on $\{T\}\times\Sc$. Fix some $(z,\theta)$ $\in$ $\Sc$, and consider some sequence 
$(t_{k},z_{k},\theta_{k})_{k\geq1}$ in  $[0,T)\times \Sc$, converging to $(T,z,\theta)$ and such that $\lim_k v(t_k,z_k,\theta_k)$ $=$ 
$v_*(T,z,\theta)$. By taking the no impulse control strategy on $[t_k,T]$, we have
\beqs
v(t_k,z_k,\theta_k) & \geq & \E\big[ U_L(Z_T^{0,t_k,z_k},\Theta_T^{0,t_k,\theta_k})\big].
\enqs
Since $(Z_T^{0,t_k,z_k},\Theta_T^{0,t_k,z_k})$ converges a.s. to $(z,\theta)$ when $k$ goes to infinity by continuity of $(Z^{0,t,z},\Theta^{0,t,\theta})$ 
in its initial condition, we deduce by Fatou's lemma that
\beq \label{vgeqU}
v_*(T,z,\theta) & \geq & U_L(z,\theta). 
\enq
On the other hand, we know from the dynamic programming QVI  that $v$ $\geq$ $\Hc v$ on $[0,T)\times\Sc$, and thus 
\beqs
v(t_k,z_k,\theta_k) &\geq&  \Hc v(t_k,z_k,\theta_k) \; \geq \; \Hc v_*(t_k,z_k,\theta_k), \;\;\; \forall k \geq 1.
\enqs 
Recalling that $\Hc v_*$ is lsc, we obtain by sending $k$ to infnity: 
\beqs
v_*(T,z,\theta) & \geq & \Hc v_*(T,z,\theta). 
\enqs
Together with \reff{vgeqU}, this proves the required viscosity supersolution property of \reff{condterm}.

We now  prove the viscosity subsolution property on $\{T\}\times\bar\Sc$, and  argue by contradiction by assuming that there exists $(\bar z,\bar \theta)$ $\in$ $\bar \Sc $ such that 
\beq \label{vsous}
\min\big[ v^*(T,\bar z,\bar \theta)  - U_L(\bar z,\bar \theta)~,~v^*(T,\bar z,\bar \theta) - \Hc v^*(T,\bar z,\bar \theta) \big] &:=& 2\eps~>~0.
\enq
One can  find a sequence of smooth functions $(\varphi^n)_{n\geq 0}$ on $[0,T]\times\bar\Sc$ such that $\varphi^n$ converges pointwisely to $v^*$ on $[0,T]\times\bar\Sc$  as $n\rightarrow\infty$. Moreover,  by \reff{vsous}  and recalling that $\Hc v^*$ is usc, we may assume that the inequality 
\beq\label{maj1}
\min\big[ \varphi^n  - U_L~,~\varphi^n - \Hc v^*  \big] &\geq& \eps,
\enq
holds on some bounded neighborhood $B^{n}$ of $(T,\bar z,\bar \theta)$ in $[0,T]\times\bar \Sc$, for $n$ large enough.  
Let  $(t_{k},z_{k},\theta_{k})_{k\geq1}$ be a sequence in $[0,T)\times \Sc$ converging to $(T,\bar z,\bar\theta)$ 
and such that $\lim_k v(t_k,z_k,\theta_k)$ $=$  $v^*(T,\bar z,\bar\theta)$.  There exists $\delta^n$ $>$ $0$ such that  
$B^{n}_k$ $:=$ $[t_{k},T]\times B(z_{k},\delta^{n})\times\Big( (\theta_{k}-\delta^{n},\theta_{k}+\delta^{n})\cap[0,T]\Big)$ $\subset$ $B^n$ for all $k$ large enough, so that \reff{maj1} holds on  $B^{n}_k$.  Since $v$ is locally bounded, there exists  some $\eta>0$ such that $|v^*|$ $\leq$ $\eta$ on $B^{n}$. We can then assume that $\varphi^n$ $\geq$ $-2\eta$ on $B^{n}$. Let us define the smooth function $\tilde \varphi _{k}^n$ on $[0,T)\times\Sc$ by
\beqs
\tilde \varphi _{k}^n(t,z,\theta) & := &  \varphi ^n(t,z,\theta)+ 4\eta\frac{|z-z_{k}|^2}{|\delta^{n}|^2}+ \sqrt{T-t}
\enqs
and observe that 
\beq\label{maj2}
(v^*-\tilde \varphi _{k}^n)(t,z,\theta) & \leq & -\eta, 
\enq
for   $ (t,z,\theta)\in[t_{k},T]\times\partial B(z_{k},\delta^{n})\times\Big( (\theta_{k}-\delta^{n},\theta_{k}+\delta^{n})\cap[0,T]\Big)$. 
Since $\Dt{\sqrt{T-t}}\longrightarrow-\infty$ as $t\rightarrow T$, we have for $k$ large enough
\beq\label{maj3}
-\Dt{\tilde \varphi _{k}^n} -\Dth{\tilde \varphi _{k}^n} - \Lc \tilde \varphi _{k}^n & \geq & 0, ~ \mbox{ on }~ B^{n}_k. 
\enq
Let $\alpha^k=(\tau^k_{j},\zeta^k_{j})_{j\geq 1}$ be a $\frac{1}{k}-$optimal control for $v(t_{k},z_{k},\theta_{k})$ with corresponding state process 
$(Z^k,\Theta^k)$, and denote by   $\sigma_n^k$ $=$  $\inf\{ s\geq t_{k}~:~ (Z_{s}^k,\Theta_{s}^k)\notin B^{n}_k \}\wedge \tau_{1}^k\wedge T$. 
From the DPP  \reff{DPP2}, this means that
\beq
v(t_{k},z_{k},\theta_{k}) - \frac{1}{k} & \leq &   \E\Big[\mathbf{1}_{\sigma_n^k<(\tau_{1}^k\wedge T)} \; 
v(\sigma_n^k,Z^k_{\sigma_n^k})\Big]  +  
\E\Big[\mathbf{1}_{\sigma_n^k=T < \tau_1^k} \;  U_{L}(Z^k_{\sigma_n^k},\Theta^k_{\sigma_n^k})\Big]
\nonumber \\
& & \; + \;  \E\Big[\mathbf{1}_{\tau_{1}^k \leq \sigma_n^k} \; v\big(\tau_1^k,\Gamma(Z_{(\tau_1^k)^-}^k,\Theta^k_{(\tau_{1}^k)^-},\zeta^k_{1}),0\big)\Big] 
\label{DPPinter}
\enq
Now,  by applying It\^{o}'s Lemma to $\tilde \varphi _{n}^k(s,Z_{s}^k,\Theta_{s}^k)$ between $t_{k}$ and $\sigma_n^k$,  
we get   from  \reff{maj1}-\reff{maj2}-\reff{maj3}, 
\beqs \label{inegphi}
\tilde \varphi^n_{k}(t_{k},z_{k},\theta_{k}) & \geq & 
\E\Big[ \mathbf{1}_{\sigma_n^k<\tau_{1}^k} \; \tilde \varphi^n_{k}(\sigma_n^k,Z^k_{\sigma_n^k},\Theta^k_{\sigma_n^k}) \Big] 
+ \E\Big[ \mathbf{1}_{\tau_{1}^k\leq \sigma_n^k}\tilde \varphi^n_{k}\big(\tau_{1}^k,Z^k_{(\tau_{1}^k)^-},\Theta^k_{(\tau_{1}^k)^-}\big) \Big]\nonumber\\
 & \geq & \E\Big[\mathbf{1}_{\sigma_n^k<(\tau_{1}^k\wedge T)}\Big( v^*(\sigma_n^k,Z^k_{\sigma_n^k},\Theta^k_{\sigma_n^k})+\eta \Big)\Big]
 +\E\Big[\mathbf{1}_{\sigma_n^k=T < \tau_1^k}\Big( U_{L}(Z^k_{\sigma_n^k},\Theta^k_{\sigma_n^k})+\eps \Big)\Big]\nonumber\\
 &  & \; + \;  \E\Big[\mathbf{1}_{\tau_{1}^k \leq \sigma_n^k} 
 \Big(v^*\big(\tau_1^k,\Gamma(Z_{(\tau_1^k)^-}^k,\Theta^k_{(\tau_{1}^k)^-},\zeta^k_{1}),0\big)  +\eps \Big) \Big]. 
\enqs
Together with \reff{DPPinter}, this implies
\beqs
\tilde \varphi^n_{k}(t_{k},z_{k},\theta_{k}) & \geq & v(t_{k},z_{k},\theta_{k}) -\frac{1}{k} +\eps\wedge\eta.
\enqs
Sending $k$, and then $n$ to infinity,  we get the required  contradiction: $v^*(T,\bar z,\bar\theta)$ $\geq$ $v^*(T,\bar z,\bar\theta) + \eps\wedge\eta$.  
\ep

\vspace{3mm}

\begin{Remark}
{\rm In order to have a complete characterization of the value function through its HJB equation, we need a uniqueness result, thus  a comparison 
principle for the QVI \reff{QVIv}-\reff{condterm}.  A key argument  originally due to Ishii \cite{ish93} for  getting  a uniqueness result for variational inequalities with impulse parts,  is to produce a strict viscosity supersolution. However, in our model, this is not possible.  Indeed, suppose  we can find a strict viscosity lsc supersolution $w$ to \reff{QVIv}, so that $(w - \Hc w)(t,z,\theta)$ 
$>$ $0$ on $[0,T)\times\Sc$.  But for $z$ $=$ $(x,y,p)$ and 
$\theta$  $=$ $0$, we have $\Gamma(z,0,e)$ $=$ $(x,y+e,p)$ for any $e$ $\Cc(z,0)$. Since  $0$ $\in$ $\Cc(z,0)$ we have $\Hc w(t,z,0)$ $=$ $\sup_{e\in [-y,0]} w(t,x,y+e,p,0)$ $\geq$ $w(t,z,0)$ $>$ $\Hc w (t,z,0)$, a contradiction.  
Actually, the main reason why one cannot obtain a strict supersolution is the absence of  fixed cost in the impulse function $\Gamma$ or in the objective functional.  
} 
\end{Remark}

\section{An approximating problem  with fixed transaction fee}

\setcounter{equation}{0} \setcounter{Assumption}{0}
\setcounter{Theorem}{0} \setcounter{Proposition}{0}
\setcounter{Corollary}{0} \setcounter{Lemma}{0}
\setcounter{Definition}{0} \setcounter{Remark}{0}

In this section, we consider a small variation of  our original model by  adding a fixed transaction fee $\eps$ $>$ $0$ 
at each trading. This means that given 
a trading strategy $\alpha$ $=$ $(\tau_n,\zeta_n)_{n\geq 0}$,  the controlled state process  $(Z=(X,Y,P),\Theta)$ jumps now at time 
$\tau_{n+1}$,  by: 
\beq \label{jumpeps}
(Z_{\tau_{n+1}},\Theta_{\tau_{n+1}}) &=& \Big( \Gamma_\eps(Z_{\tau_{n+1}^-},\Theta_{\tau_{n+1}^-},\zeta_{n+1}), 0 \Big), 
\enq
where $\Gamma_\eps$ is the function defined on $\R\times\R_+\times\R_+^*\times [0,T]\times\R$ into 
$\R\cup\{-\infty\}\times\R\times\R_+^*$ by:
\beqs
\Gamma_\eps(z,\theta,e) &=& \Gamma(z,\theta,e) - (\eps,0,0) \; = \; 
\Big( x - e p f(e,\theta) - \eps, y+ e, p\Big),
\enqs
for $z$ $=$ $(x,y,p)$. The dynamics of $(Z,\Theta)$ between trading dates is given as before. 
We also introduce a modified liquidation function $L_\eps$ defined  by: 
\beqs
L_{\eps}(z,\theta) & = & \max[ x ,L(z,\theta) - \eps ], \;\;\; (z,\theta) = (x,y,p,\theta) \in \R\times\R_+\times\R_+^*\times [0,T]. 
\enqs 
The interpretation of this modified liquidation function is the following. Due to the presence of the transaction fee at each trading,  it may be advantageous for the investor not to liquidate his position in stock shares (which would give him $L(z,\theta)-\eps$), and rather bin his stock shares, 
by keeping only his cash amount (which would give him $x$).  Hence, the investor chooses 
 the best of these two possibilities, which induces a liquidation  value $L_\eps(z,\theta)$. 

We then introduce  the corresponding   solvency region $\Sc_\eps$ with its closure $\bar\Sc_\eps$ $=$ $\Sc_\eps$ $\cup$ $\partial\Sc_\eps$, and boundary 
$\partial\Sc_\eps$ $=$ $\partial_y\Sc_\eps$ $\cup$ $\partial_L\Sc_\eps$:
\beqs
\Sc_{\eps} & = & \Big\{  (z,\theta) =(x,y,p,\theta)\in\R\times\R_{+}\times\R_{+}^*\times [0,T]:~ y  > 0 \; \mbox{ and } L_{\eps}(z,\theta)>0  \Big\}, \\
\partial_y\Sc_\eps &=& \Big\{  (z,\theta) =(x,y,p,\theta)\in\R\times\R_{+}\times\R_{+}^*\times [0,T]:~ y  = 0 \; \mbox{ and } 
\; L_{\eps}(z,\theta) \geq 0  \Big\}, \\
\partial_L\Sc_\eps &=& \Big\{  (z,\theta) =(x,y,p,\theta)\in\R\times\R_{+}\times\R_{+}^*\times\R_{+}:~ L_{\eps}(z,\theta)=0  \Big\}. 
\enqs
We also introduce  the corner lines of $\partial\Sc_\eps$. For simplicity of presentation, we consider a temporary price impact function $f$ in the form: 
\beqs
f(e,\theta) &=& \tilde f \Big( \frac{e}{\theta} \Big) ~=~\exp\Big(\lambda \frac{e}{\theta}\Big)  
\Big( \kappa_{a}\mathbf{1}_{e>0} +\mathbf{1}_{e=0}  +\kappa_{b}\mathbf{1}_{e<0} \Big)\mathbf{1}_{\theta>0},
\enqs
where $0<\kappa_b<1<\kappa_a$, and $\lambda$ $>$ $0$. A straightforward analysis of the function $L$ shows that $y$ $\mapsto$ 
$L(x,y,p,\theta)$ is increasing on $[0,\theta/\lambda]$, decreasing on $[\theta/\lambda,\infty)$ with $L(x,0,p,\theta)$ $=$ $x$ $=$ 
$L(x,\infty,p,\theta)$, and $\max_{y>0} L(x,y,p,\theta)$ $=$ $L(x,\theta/\lambda,p,\theta)$ $=$ $x+p\frac{\theta}{\lambda}\tilde f(-1/\lambda)$. 
We first get the form of the sets $\Cc(z,\theta)$:
\beqs
\Cc(z,\theta) & = & [-y,\bar e(z,\theta)]\;,
\enqs
where the function $\bar e$ is defined in Lemma \ref{lemsetadmi}.
 We then distinguish two cases: 
(i) If $p\frac{\theta}{\lambda}\tilde f(-1/\lambda)$ $<$ $\eps$, then $L_\eps(x,y,p,\theta)$ $=$ $x$. (ii) If 
$p\frac{\theta}{\lambda}\tilde f(-1/\lambda)$ $\geq$ $\eps$, then there exists an unique $y_1(p,\theta)$ $\in$ $(0,\theta/\lambda]$ and 
$y_2(p,\theta)$ $\in$ $[\theta/\lambda,\infty)$ such that $L(x,y_1(p,\theta),p,\theta)$ $=$ $L(x,y_2(p,\theta),p,\theta)$ $=$ $x$, and 
$L_\eps(x,y,p,\theta)$ $=$ $x$ for $y$ $\in$ $[0,y_1(p,\theta))\cup (y_2(p,\theta),\infty)$, $L_\eps(x,y,p,\theta)$ $=$ $L(x,y,p,\theta)-\eps$ for 
$y$ $\in$ $[y_1(p,\theta),y_2(p,\theta)]$.  We then denote by 
\beqs
D_0  & = & \{0\}\times\{0\}\times\R_+^*\times [0,T] \; = \;  \partial_y\Sc_\eps\cap\partial_L\Sc_\eps,  \\
D_{1,\eps} &=&  \Big\{ (0,y_1(p,\theta),p,\theta):  p\frac{\theta}{\lambda}\tilde f\big(\frac{-1}{\lambda}\big) \geq \eps, \; \theta \in [0,T] \Big\}, \\
D_{2,\eps} &=& \Big\{ (0,y_2(p,\theta),p,\theta):  p\frac{\theta}{\lambda} \tilde f\big(\frac{-1}{\lambda}\big) \geq \eps, \; \theta \in [0,T] \Big\}. 
\enqs
Notice that  the inner normal vectors at the corner lines $D_{1,\eps}$ and $D_{2,\eps}$  form an acute angle (positive scalar product), while 
we have a right angle at the corner $D_0$.

\begin{figure}
\begin{center}\label{Seps2d}
\includegraphics[angle= -90,width=16cm]
{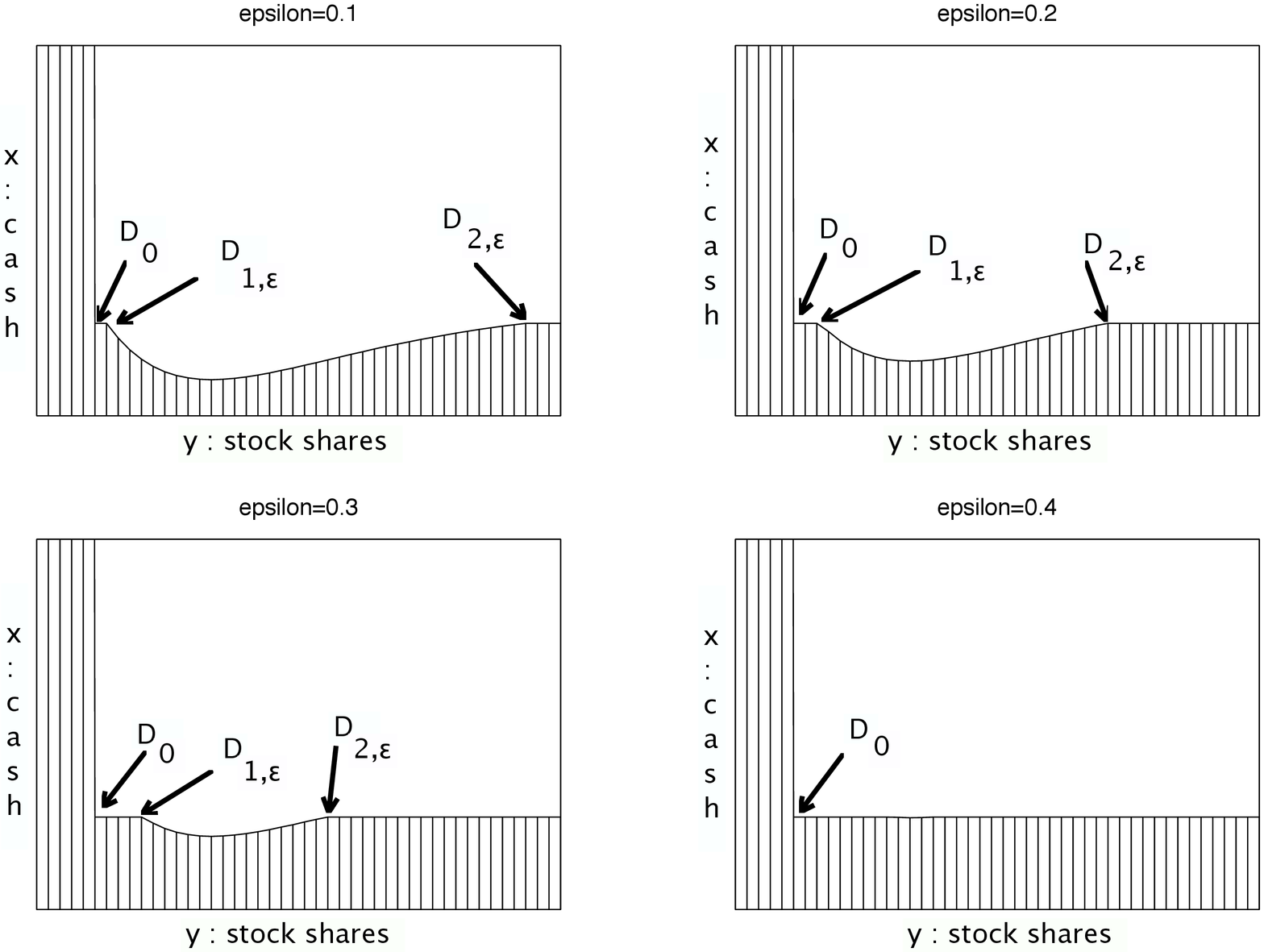}
\caption{Domain $\Sc_{\eps}$ in the nonhatched zone for fixed $p=1$ and $\theta=1$ and $\eps$ evolving from $0.1$ to $0.4$. Here $\kappa_{b}=0.9$ and $f(e,\theta)=\kappa_{b}\exp\big(\frac{e}{\theta}\big)$ for $e<0$. Notice that for $\eps$ large enough, $\Sc_{\eps}$ is equal to open orthant $\R_{+}^*\times\R_{+}^*$.
}
\end{center}

\end{figure}

\begin{figure}
\label{Seps3d-theta-fix}
\begin{center}
\includegraphics[width=16cm,height=9cm]{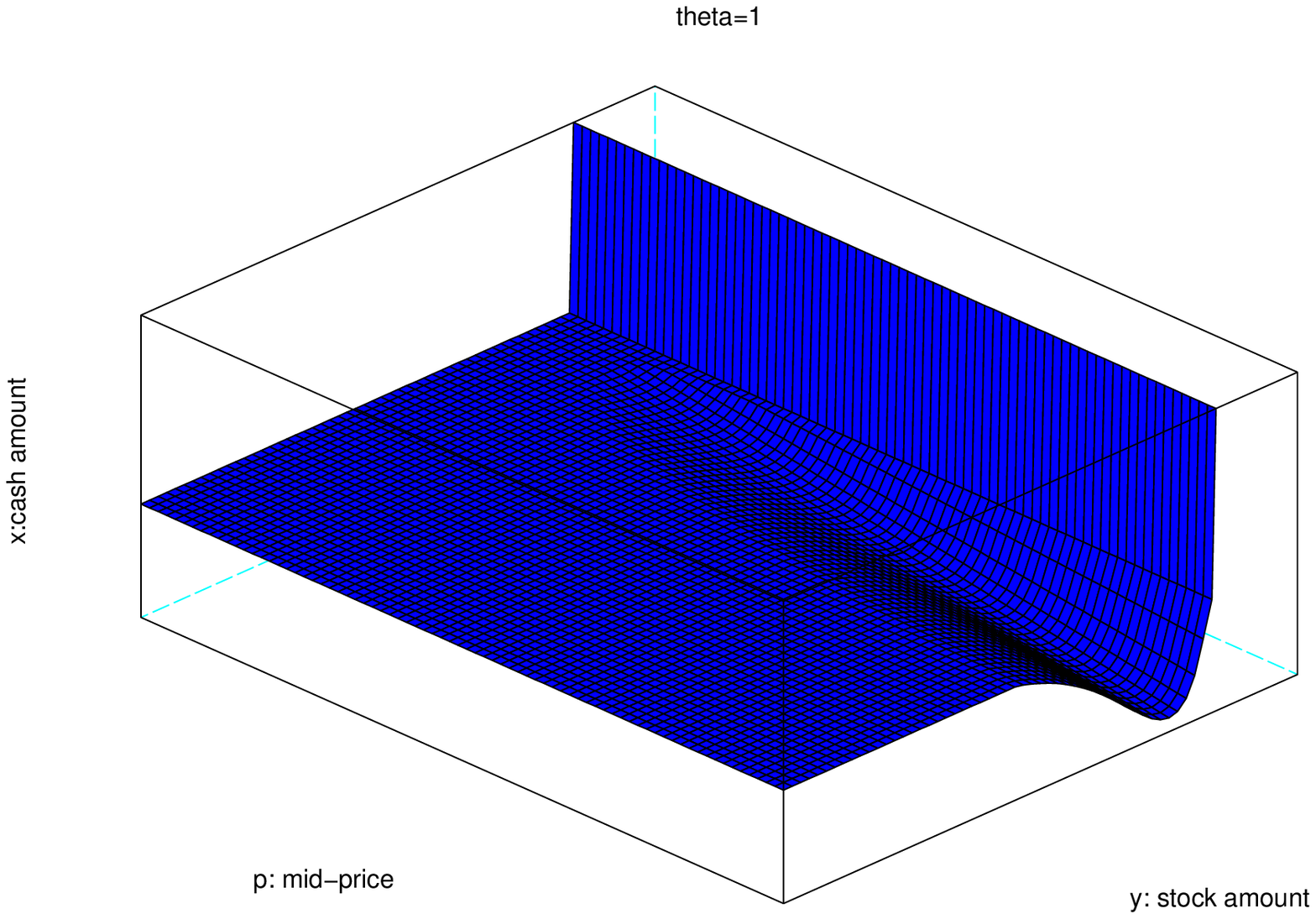}
\caption{Lower bound of the domain $\Sc_{\eps}$ for fixed  $\theta=1$ and $f(e,\theta)=\kappa_{b}\exp\big(\frac{e}{\theta}\big)$ for $e<0$. Notice that when $p$ is fixed, we obtain the Figure 4.
}
\end{center}
\end{figure}

\begin{figure}
\begin{center}\label{Seps3d-p-fix}
\includegraphics[width=16cm,height=9cm]{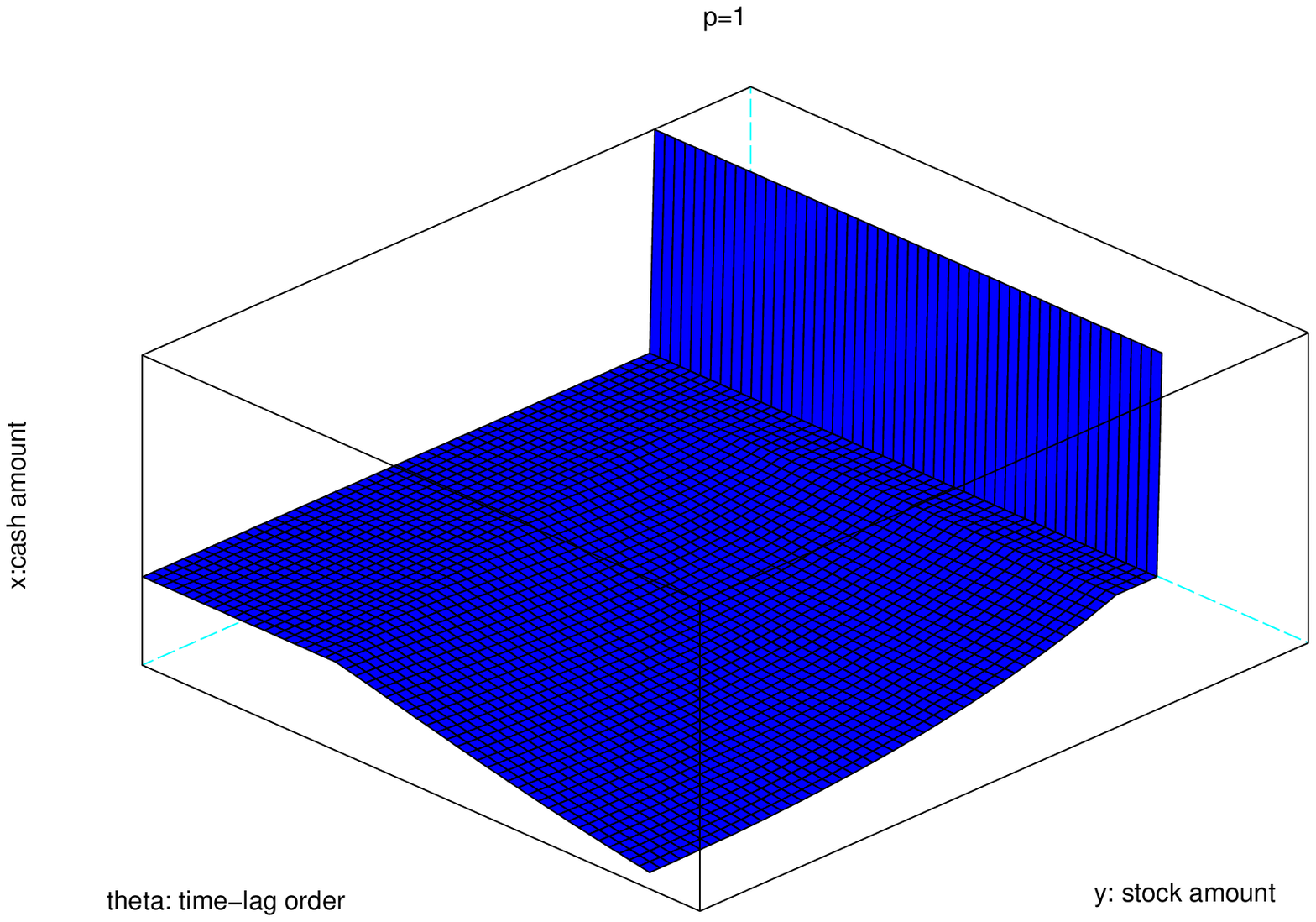}
\caption{Lower bound of the domain $\Sc_{\eps}$ for fixed  $p=1$ and $\eps=0.2$. Here $\kappa_{b}=0.9$ and $f(e,\theta)=\kappa_{b}\exp\big(\frac{e}{\theta}\big)$ for $e<0$. Notice that when $\theta$ is fixed, we obtain the Figure 4. 
}
\end{center}
\end{figure}

\vspace{1mm}

Next,   we define  the set of admissible trading strategies as follows. 
Given $(t,z,\theta)$ $\in$ $[0,T]\times\bar\Sc_\eps$, we say that the impulse control  $\alpha$ is admissible, denoted by $\alpha$ $\in$ 
$\Ac^\eps(t,z,\theta)$, if $\tau_0$ $=$ $t-\theta$, $\tau_n$ $\geq$ $t$, $n$ $\geq$ $1$,  and the controlled state process $(Z^\eps,\Theta)$  
solution to  \reff{dynTheta}-\reff{Y}-\reff{dynP}-\reff{Xconst}-\reff{jumpeps}, with  an initial state 
$(Z_{t^-}^\eps,\Theta_{t^-})$ $=$ $(z,\theta)$ (and the convention that  $(Z_{t}^\eps,\Theta_{t})$ $=$ $(z,\theta)$ if $\tau_1$ $>$ $t$),      
satisfies  $(Z_{s}^\eps,\Theta_s)$ $\in$ $[0,T]\times\bar\Sc_\eps$ for all $s$ $\in$ $[t,T]$.  Here, we stress the dependence of 
$Z^\eps$ $=$ $(X^\eps,Y,P)$ in  $\eps$ appearing in the transaction function  $\Gamma_\eps$, and we notice that it affects only the cash component. 
Notice that $\Ac^\eps(t,z,\theta)$ is nonempty for any 
$(t,z,\theta)$ $\in$ $[0,T]\times\bar\Sc_\eps$. Indeed, for $(z=(x,y,p),\theta)$ $\in$ $\bar\Sc_\eps$, i.e. $L_\eps(z,\theta)$ $=$ 
$\max(x,L(z,\theta)-\eps)$ $\geq$ $0$, we distinguish two cases:  (i)  if $x$ $\geq$ $0$, then by doing none transaction, the associated state process 
$(Z^\eps=(X^\eps,Y,P),\Theta)$ satisfies $X_s^\eps$ $=$ $x$ $\geq$ $0$, $t\leq s\leq T$, and thus this  zero transaction is admissible; 
(ii) if $L(z,\theta)-\eps$ $\geq$ $0$, then by liquidating immediately all the stock shares, and doing nothing more after, the associated state process 
satisfies $X_s^\eps$ $=$ $L(z,\theta)-\eps$, $Y_s$ $=$ $0$,   and thus $L_\eps(Z_s^\eps,\Theta_s)$ $=$ $X_s^\eps$ $\geq$ $0$, $t\leq s\leq T$, which shows that this immediate transaction is  admissible.

Given the utility function $U$ on $\R_+$, and the liquidation utility function defined on $\bar\Sc_\eps$ by 
$U_{L_\eps}(z,\theta)$ $=$ $U(L_\eps(z,\theta))$, 
we then consider  the associated  optimal portfolio liquidation problem defined via its value function by:
 \beq \label{defveps}
 v_{\varepsilon}(t,z,\theta) & = & \sup_{\alpha\in\Ac^\varepsilon(t,z,\theta)}\E\big[ U_{L_{\eps}}(Z_{T}^\eps,\Theta_T) \big], \;\;\; 
 (t,z,\theta) \in [0,T]\times \bar\Sc_\eps.
 \enq

Notice that  when $\eps$ $=$ $0$,  the above problem reduces to the optimal portfolio liquidation problem described in Section 2, and 
in particular $v_0$ $=$ $v$. 
The main purpose of this section is to provide a unique PDE characterization of the value functions $v_\eps$, $\eps$ $>$ $0$, and to prove that 
 the sequence $(v_\eps)_\eps$ converges to the original value function $v$ as $\eps$ goes to zero.

 \vspace{2mm}

We define the set of admissible transactions in the model with fixed transaction fee by:
 \beqs
 \Cc_{\eps}(z,\theta) & = & \Big\{e\in\R: \Big(\Gamma_\eps(z,\theta,e),0\Big) \in \bar\Sc_\eps  \Big\}, \;\;\; (z,\theta) \in \bar\Sc_\eps.
 \enqs 
A similar calculation as in Lemma \ref{lemsetadmi} shows that  for $(z=(x,y,p),\theta)$ $\in$ $\bar\Sc_\eps$, 
\beqs
\Cc_\eps(z,\theta) &=& \left\{
				    \begin{array}{cl}
				    [-y,\bar e_\eps(z,\theta)], & \mbox{ if } \theta > 0 \mbox{ or } x \geq \eps, \\
				     \emptyset, & \mbox{ if } \; \theta = 0 \; \mbox{ and } x < \eps,
				    \end{array}
				    \right.	
\enqs
where  $\bar e(z,\theta)$ $=$ $\sup\{ e \in \R: ep\tilde f(e/\theta) \leq x- \eps\}$ if $\theta$ $>$ $0$ and $\bar e(z,0)$ $=$ $0$ if $x$ $\geq$ $\eps$. 
Here,  the set $[-y,\bar e_\eps(z,\theta)]$ should be viewed as empty when $\bar e(z,\theta)$ $<$ $y$, i.e. 
$x+py\tilde f(-y/\theta)-\eps$ $<$ $0$. 
We also easily check that $\Cc_\eps$ is continuous for the Hausdorff metric. We then consider the impulse operator $\Hc_\eps$ by 
\beqs
\Hc_\eps w(t,z,\theta) &=& \sup_{e \in \Cc_\eps(z,\theta)} w(t,\Gamma_\eps(z,\theta,e),0), \;\;\; (t,z,\theta) \in [0,T]\times\bar\Sc_\eps,
\enqs
for any locally bounded function $w$ on $[0,T]\times\bar\Sc_\eps$, with the convention that $\Hc_\eps w(t,z,\theta)$ $=$ $-\infty$ when 
$\Cc_\eps(z,\theta)$ $=$ $\emptyset$.

Next, consider again the Merton liquidation  function $L_M$, and observe similarly as in \reff{LMgamexpli} that
\beq
L_M(\Gamma_\eps(z,\theta,e)) - L_M(z) &=&  e p \Big( 1 - f\big(e,\theta\big) \Big)  - \eps  \nonumber \\
& \leq & - \eps, \;\;\;\;\;  \forall (z,\theta)\in\bar\Sc_\eps, \;\; e \in \R.  \label{LMeps}
\enq
This implies in particular that 
\beq \label{HcLM}
\Hc_\eps L_M & < & L_M \;\;\; \mbox{ on } \;  \bar\Sc_\eps. 
\enq
Since $L_\eps$ $\leq$ $L_M$,  we observe from \reff{LMeps} 
that if  $(z,\theta)$ $\in$ $\Nc_\eps$ $:=$ $\{ (z,\theta) \in \bar\Sc_\eps: L_M(z) < \eps\}$, then 
$\Cc_\eps(z,\theta)$ $=$ $\emptyset$. Moreover, we deduce from \reff{LMeps} that for all $\alpha$ $=$ $(\tau_n,\zeta_n)_{n\geq 0}$ 
$\in$ $\Ac^\eps(t,z,\theta)$ associated to the state process $(Z,\Theta)$, $(t,z,\theta)$ $\in$ $[0,T]\times\bar\Sc_\eps$: 
\beqs
0 \; \leq \; L_M(Z_T) &=& L_M(Z_T^{0,t,z}) + \sum_{n\geq 0}  \Delta L_M(Z_{\tau_n}) \\
& \leq &  L_M(Z_T^{0,t,z}) - \eps N_T(\alpha),
\enqs
where we recall that $N_T(\alpha)$ is the number of trading times over the whole horizon $T$.  This shows that 
\beqs
N_T(\alpha) & \leq & \frac{1}{\eps}  L_M(Z_T^{0,t,z})\; < \; \infty \;\;\; a.s. 
\enqs
In other words,  we  see  that, under the presence of fixed transaction fee,  
the number of intervention  times over a finite interval for an admissible  trading strategy   is finite almost surely.

\vspace{2mm}

The dynamic programming equation associated to the control problem \reff{defveps} is
 \beq
 \min\Big[ -\Dt{w} - \Dth{w} - \Lc w~,~w -\Hc_\eps w \Big] & = & 0, \;\;\;\; \mbox{ in } [0,T)\times\bar\Sc_\eps, \label{QVIvareps} \\
 \min\big[ w -  U_{L_\eps}~,~w - \Hc_\eps w \big] &=& 0,  \;\;\;\; \mbox{ in }  \{T\} \times\bar\Sc_\eps. \label{termeps}
\enq

The main result of this section is stated as follows.

 \begin{Theorem} \label{theomainapprox}
(1) The sequence $(v_\eps)_\eps$ is nonincreasing, and converges  pointwise on 
$[0,T]\times(\bar\Sc\setminus\partial_L\Sc)$ towards  $v$  as $\eps$ goes to zero. 

\noindent (2)  For any  $\eps$ $>$ $0$, the value function $v_\eps$ is continuous on $[0,T)\times\Sc_\eps$, and is  the unique 
(in $[0,T)\times\Sc_\eps)$ constrained viscosity solution to \reff{QVIvareps}-\reff{termeps}, satisfying the growth condition:
\beq \label{growthveps}
|v_\eps(t,z,\theta)| & \leq & K (1+ L_M(z)^\gamma), \;\;\; \forall (t,z,\theta) \in [0,T]\times\bar\Sc_\eps,
\enq
for some positive constant $K$, and the boundary condition: 
\beq
\lim_{(t',z',\theta')\rightarrow (t,z,\theta)} v_\eps(t',z',\theta') &=& v(t,z,\theta) \nonumber \\
& = & U(0), \;\;\;  \forall (t,z=(0,0,p),\theta) \in [0,T]\times D_0.  \label{bounveps}
\enq
\end{Theorem}

\vspace{2mm}

We  first prove the convergence of the sequence of value functions $(v_\eps)$.

\vspace{2mm}

\noindent  {\bf Proof of Theorem \ref{theomainapprox} (1)}.

\noindent   Notice that for any $0<\eps_1\leq\eps_2$, we have 
$L_{\eps_2}$ $\leq$ $L_{\eps_1}$ $\leq$ $L$,   $\Ac^{\eps_2}(t,z,\theta)$ $\subset$ 
 $\Ac^{\eps_1}(t,z,\theta)$ $\subset$ $\Ac(t,z,\theta)$,   for $t$ $\in$ $[0,T]$, $(z,\theta)$ $\in$ $\bar\Sc_{\eps_2}$ $\subset$ $\bar\Sc_{\eps_1}$ 
 $\subset$ $\bar\Sc$,  and  for $\alpha$ $\in$ $\Ac^{\eps_2}(t,z,\theta)$, 
  $L_{\eps_2}(Z^{\eps_2},\Theta)$ $\leq$  $L_{\eps_2}(Z^{\eps_1},\Theta)$ $\leq$ $L_{\eps_1}(Z^{\eps_1},\Theta)$ $\leq$ $L(Z,\Theta)$.  This shows that  the sequence $(v_\eps)$ is nonincreasing, and is upper-bounded by the  value function $v$ without transaction fee, so that 
\beq \label{liminfveps}
\lim_{\eps\downarrow 0} v_\eps (t,z,\theta) &\leq& v(t,z,\theta), \;\;\; \forall (t,z,\theta) \in [0,T]\times\bar\Sc.  
\enq

Fix now some point $(t,z,\theta)$ $\in$ $[0,T]\times(\bar\Sc\setminus\partial_L\Sc)$.  From the representation \reff{vfini2} of $v(t,z,\theta)$, 
there exists for any  $n$ $\geq$ $1$, an $1/n$-optimal control $\alpha^{(n)}$ $=$ $(\tau_k^{(n)},\zeta_k^{(n)})_k$ 
$\in$ $\Ac_{\ell_+}^b(t,z,\theta)$ with asso\-ciated state process 
$(Z^{(n)}=(X^{(n)},Y^{(n)},P),\Theta^{(n)})$ and number of trading times $N^{(n)}$: 
\beq \label{voptin}
\E\big[ U(X_{T}^{(n)}) \big] & \geq & v(t,z,\theta) - \frac{1}{n}.
\enq
We denote by  $(Z^{\eps,(n)},\Theta^{(n)})$ $=$ $(X^{\eps,(n)},Y^{(n)},P),\Theta^{(n)})$ the state process controlled by $\alpha^{(n)}$ in the model 
with transaction fee $\eps$ (only the cash component is affected by $\eps$), and we observe that for all $t\leq s\leq T$, 
\beq \label{XNeps}
X_s^{\eps,(n)} &=& X_s^{(n)} - \eps N_s^{(n)} \; \nearrow \; X_s^{(n)},  \;\;\; \mbox{ as } \; \eps \; \mbox{ goes to zero.}  
\enq
Given $n$, we  consider the family of  stopping times:
\beqs
\sigma^{(n)}_\eps & = & \inf\big\{ s\geq t:~ L(Z_{s}^{\varepsilon,(n)},\Theta_s^{(n)}) \leq \varepsilon \big\} \wedge T, \;\;\;\;\; \eps > 0. 
\enqs 
Let us prove that 
\beq \label{limsig}
\lim_{\eps\searrow 0} \sigma^{(n)}_\eps & = & T  \;\; a.s. 
\enq
Observe that for $0<\eps_1\leq\eps_2$, $X_s^{\eps_2,(n)}$ $\leq$ $X_s^{\eps_1,(n)}$, and so 
$L(Z_s^{\eps_2,(n)},\Theta_s)$ $\leq$  $L(Z_s^{\eps_1,(n)},\Theta_s)$ for $t\leq s\leq T$.  This implies clearly that the sequence 
$(\sigma_\eps^{(n)})_\eps$ is nonincreasing. Since this sequence is bounded by $T$, it admits a limit, denoted by $\sigma_0^{(n)}$ $=$ 
$\lim_{\eps\downarrow 0} \uparrow \sigma_\eps^{(n)}$. 
Now, by definition of $\sigma^{(n)}_\eps$, we have  $L(Z_{\sigma^{(n)}_\eps}^{\varepsilon,(n)},\Theta_{\sigma^{(n)}_\eps}^{(n)})$ $\leq$ 
$\varepsilon$,  for all $\eps$ $>$ $0$.  By sending $\eps$ to zero, we then get with \reff{XNeps}:
\beqs
L(Z_{\sigma^{(n),-}_0}^{(n)},\Theta_{\sigma^{(n),-}_0}^{(n)}) & = &  0 \;\;\; a.s. 
\enqs
Recalling the definition  of  $\Ac_{\ell_+}^b(t,z,\theta)$, this implies that 
$\sigma^{(n)}_0$ $=$ $\tau_{k}^{(n)}$ for some $k\in\{1,\ldots,N^{(n)}+1\}$ with the convention $\tau_{N^{(n)}+1}^{(n)}$ $=$ $T$. 
If $k$ $\leq$ $N^{(n)}$, arguing as in \reff{Lalpha}, we get a contradiction with the solvency constraints. Hence we get $\sigma^{(n)}_0$ $=$ $T$.  

Consider now the trading strategy 
$\tilde\alpha^{\eps,(n)}$ $\in$  $\Ac$ consisting in following $\alpha^{(n)}$ until time $\sigma^{(n)}_\eps$ and liquidating all the stock shares at time $\sigma_\eps^{(n)}$, i.e. 
\beqs
\tilde\alpha^{\eps,(n)} &=& (\tau_k^{(n)},\zeta_k^{(n)}) \mathbf{1}_{\tau_k < \sigma_\eps^{(n)}} \cup (\sigma_\eps^{(n)},-Y_{\sigma_\eps^{(n),-}}). 
\enqs
We denote by $(\tilde Z^{\eps,(n)}=(\tilde X^{\eps,(n)},\tilde Y^{\eps,(n)},P),\tilde\Theta^{\eps,(n)})$ the associated state process in the market with transaction fee $\eps$. 
By construction, we have for all 
$t$ $\leq$ $s$ $<$ $\sigma^{(n)}_\eps$:  $L(\tilde Z_s^{\eps,(n)},\tilde\Theta_s^{\eps,(n)})$ $=$ $L(Z_s^{\eps,(n)},\Theta_s^{(n)})$ 
$\geq$ $\eps$, and thus $L_\eps(\tilde Z_s^{\eps,(n)},\tilde\Theta_s^{\eps,(n)})$ $\geq$ $0$. 
At the transaction time $\sigma^{(n)}_\eps$, we then have 
$\tilde X^{\eps,(n)}_{\sigma^{(n)}_\eps}$ $=$  $L(\tilde Z_{\sigma^{(n),-}_\eps}^{\eps,(n)},\tilde\Theta_{\sigma^{(n),-}_\eps}^{\eps,(n)}) -\eps$ 
$=$  $L(Z_{\sigma^{\eps,(n),-}_\eps}^{(n)},\Theta_{\sigma^{(n),-}_\eps}^{(n)}) - \eps$,  $\tilde Y^{\eps,(n)}_{\sigma^{(n)}_\eps}$ $=$ $0$. 
After time $\sigma_\eps^{(n)}$, there is no more transaction in $\tilde\alpha^{\eps,(n)}$, and so 
\beq
\tilde X_s^{\eps,(n)} \; = \;  \tilde X^{\eps,(n)}_{\sigma^{(n)}_\eps} &=& 
L(Z_{\sigma^{\eps,(n),-}_\eps}^{(n)},\Theta_{\sigma^{(n),-}_\eps}^{(n)}) - \eps \; \geq \; 0,  \label{Xintersig}  \\
\tilde Y_s^{\eps,(n)}  \; = \; \tilde Y^{\eps,(n)}_{\sigma^{(n)}_\eps} &=&  0,  \;\;\;\;\;\;\;\;  \sigma_\eps^{(n)}\leq s \leq T, \label{Yintersig}
\enq
and thus $L_\eps(\tilde Z_s^{\eps,(n)},\tilde\Theta_s^{\eps,(n)})$ $=$ $\tilde X_s^{\eps,(n)}$ $\geq$ $0$ for $\sigma_\eps^{(n)}\leq s \leq T$. 
This shows that  $\tilde\alpha^{\eps,(n)}$ lies in $\Ac^\eps(t,z,\theta)$, and thus by definition of $v_\eps$: 
\beq \label{interveps}
v_{\eps}(t,z) & \geq &  \E\big[ U_{L_{\eps}}\big(\tilde Z_{T}^{\eps,(n)},\tilde\Theta_T^{\eps,(n)}\big) \big]. 
\enq
Let us check that given $n$, 
\beq \label{limepsL}
\lim_{\eps\downarrow 0} L_\eps\big(\tilde Z_{T}^{\eps,(n)},\tilde\Theta_T^{\eps,(n)}\big) &=& X_T^{(n)}, \;\;\;\;\; a.s. 
\enq
To alleviate notations, we set $N$ $=$ $N_T^{(n)}$ the total number of trading times of $\alpha^{(n)}$. 
If the last trading time of $\alpha^{(n)}$ occurs strictly before $T$, then  we do not trade anymore until the final horizon $T$, and so 
\beq \label{intertaun}
X_T^{(n)} \; = \;  X_{\tau_{N}}^{(n)}, & \mbox{and}   & Y_T^{(n)} \; = \; Y_{\tau_N}^{(n)} \; = \; 0, \;\;\; \mbox{ on } \; \{ \tau_N < T \}.
\enq
By \reff{limsig},  we have for $\eps$ small enough: 
$\sigma_\eps^{(n)}$ $>$ $\tau_N$, and so $\tilde X_{\sigma_\eps^{(n),-}}^{\eps,(n)}$ $=$ $X_{\tau_N}^{\eps,(n)}$, 
$\tilde Y_{\sigma_\eps^{(n),-}}^{\eps,(n)}$ $=$ $Y_{\tau_N}^{(n)}$ $=$ $0$.  The final liquidation at time $\sigma_\eps^{(n)}$  yields:  
$\tilde X_T^{\eps,(n)}$ $=$ $\tilde X_{\sigma_\eps^{(n)}}^{\eps,(n)}$ $=$ $\tilde X_{\sigma_\eps^{(n),-}}^{\eps,(n)}-\eps$ $=$ 
$X_{\tau_N}^{\eps,(n)} -\eps$, and  $\tilde Y_T^{\eps,(n)}$ $=$ $\tilde Y_{\sigma_\eps^{(n)}}^{\eps,(n)}$ $=$ $0$.  We then obtain
\beqs
L_\eps \big(\tilde Z_{T}^{\eps,(n)},\tilde\Theta_T^{\eps,(n)}\big) &=& 
\max\Big( \tilde X_T^{\eps,(n)} , L\big(\tilde Z_{T}^{\eps,(n)},\tilde\Theta_T^{\eps,(n)}\big) - \eps \Big) \nonumber \\
&=& \tilde X_T^{\eps,(n)} \; = \;  X_{\tau_N}^{\eps,(n)} -\eps \;\;\;\;\;  \mbox{ on } \; \{ \tau_N < T \} \nonumber \\
&=& X_T^{(n)} - (1+N) \eps \;\;\;\;\;  \mbox{ on } \; \{ \tau_N < T \}, \label{LepsT}
\enqs
by \reff{XNeps} and \reff{intertaun}, which shows that the convergence in \reff{limepsL} holds on $\{\tau_N<T\}$. 
 If the last trading of $\alpha^{(n)}$ occurs at time $T$, this means that we liquidate all stock shares at $T$, and so 
\beq \label{Xtaun<T}
X_T^{(n)} &=& L(Z_{T^-}^{(n)},\Theta_{T^-}^{(n)}), \;\; Y_T^{(n)} \; = \; 0   \;\;\;\;\;  \mbox{ on } \; \{ \tau_N =T \}. 
\enq
On the other hand, by \reff{Xintersig}-\reff{Yintersig}, we have
\beqs
L_\eps \big(\tilde Z_{T}^{\eps,(n)},\tilde\Theta_T^{\eps,(n)}\big) \; = \;   \tilde X_T^{\eps,(n)} &  = & 
L(Z_{\sigma^{\eps,(n),-}_\eps}^{(n)},\Theta_{\sigma^{(n),-}_\eps}^{(n)}) - \eps  \\
& \longrightarrow &  L(Z_{T^-}^{(n)},\Theta_{T^-}^{(n)}) \;\;\;\;\; \mbox{ as } \eps \; \mbox{ goes to zero},
\enqs
by \reff{limsig}.   Together with \reff{Xtaun<T}, this implies  that the convergence in \reff{limepsL} also holds on $\{\tau_N=T\}$, and 
thus  almost surely.  Since $0\leq L_\eps \leq L$, we immediately see by Proposition \ref{PropBound} that the sequence 
$\{ U_{L_{\eps}}\big(\tilde Z_{T}^{\eps,(n)},\tilde\Theta_T^{\eps,(n)}\big), \eps > 0\}$  is uniformly integrable, so that by sending $\eps$ to zero in 
\reff{interveps} and using \reff{limepsL}, we get
\beqs
\lim_{\eps\downarrow 0} v_\eps(t,z,\theta) & \geq &  \E\big[ U(X_{T}^{(n)}) \big] \; \geq \; v(t,z) - \frac{1}{n},
\enqs
from \reff{voptin}. By sending $n$ to infinity, and recalling \reff{liminfveps}, this completes the proof of assertion (1) in  Theorem \ref{theomainapprox}. 
\ep

\vspace{2mm}

We now   turn to the viscosity characterization of $v_\eps$. The viscosity property of $v_\eps$ is proved similarly as for $v$, and is then omitted. 
From Proposition \ref{PropBound}, and since $0$ $\leq$ $v_\eps$ $\leq$ $v$, we know that the value functions $v_\eps$  lie in the set of functions satisfying the growth condition in \reff{growthveps}, i.e. 
\beqs
\Gc_\gamma([0,T]\times\bar\Sc_\eps) &=& \Big\{ w: [0,T]\times\bar\Sc_\eps \rightarrow \R, \; 
\sup_{[0,T]\times\bar\Sc_\eps}  \frac{|w(t,z,\theta)|}{1+ L_M(z)^\gamma} \; < \; \infty \Big\}.
\enqs 
The boundary property \reff{bounveps} is immediate. Indeed, fix $(t,z=(x,0,p),\theta)$ $\in$ $[0,T]\times\partial_y\Sc_\eps$, and consider 
an arbitrary  sequence  $(t_n,z_n=(x_n,y_n,p_n),\theta_n)_n$ in $[0,T]\times\bar\Sc_\eps$ converging to $(t,z,\theta)$. Since $0$ $\leq$ 
$L_\eps(z_n,\theta_n)$ $=$ $\max(x_n,L(z_n,\theta_n)-\eps)$, and $y_n$ goes to zero, this implies that for $n$ large enough, $x_n$ $=$ 
$L_\eps(z_n,\theta_n)$  $\geq$ $0$. 
By considering from $(t_n,z_n,\theta_n)$ 
the admissible strategy of doing none transaction, which leads to a final liquidation value $X_T$ $=$ $x_n$, we have 
$U(x_n)$ $\leq$ $v_\eps(t_n,z_n,\theta_n)$ $\leq$ $v(t_n,z_n,\theta_n)$.  
Recalling Corollary \ref{coroSy}, we then obtain  the continuity of $v_\eps$ on $\partial_y\Sc_\eps$ with 
$v_\eps(t,z,\theta)$ $=$ $U(x)$ $=$ $v(t,z,\theta)$ for $(z,\theta)$ $=$ $(x,0,p,\theta)$ $\in$ $\partial_y\Sc_\eps$, and in particular  \reff{bounveps}. 
Finally, we address the uniqueness issue, which is a direct consequence of the following 
comparison principle for constrained (discontinuous) viscosity solution to  \reff{QVIvareps}-\reff{termeps}.

\begin{Theorem}\label{Thmcompeps} (Comparison principle) 

\noindent Suppose $u$ $\in$ $\Gc_{\gamma}([0,T]\times\bar\Sc _{\eps})$ is a usc viscosity subsolution to \reff{QVIvareps}-\reff{termeps} 
on $[0,T]\times\bar\Sc_\eps$, and $w$  $\in$ $\Gc_{\gamma}([0,T]\times\bar\Sc _{\eps})$ is a lsc  viscosity supersolution to 
\reff{QVIvareps}-\reff{termeps}  on $[0,T]\times\Sc_\eps$ such that
\beq \label{uwboun}
u(t,z,\theta) & \leq & \liminf_{\tiny{\begin{array}{c}(t',z',\theta')\rightarrow(t,z,\theta)\\ (t',z',\theta')\in[0,T)\times\Sc_{\eps} \end{array}}}w(t',z',\theta'), \quad\forall(t,z,\theta)\in [0,T] \times D_0. 
\enq
Then, 
\beq \label{uleqw}
u & \leq & w \quad\mbox{ on }\quad [0,T]\times \Sc_\eps. 
\enq
\end{Theorem}

Notice that with respect to usual comparison principles for parabolic PDEs where we compare a viscosity subsolution and a viscosity supersolution from the inequalities on the domain and at the terminal date,  we require here in addition a comparison on the boundary 
$D_0$  due to the  non smoothness  of the domain $\bar\Sc_\eps$ on this right angle of the boundary. 
A similar feature appears also in \cite{lyvmnipha07}, and we shall only emphasize the main arguments adapted from \cite{bar94}, 
for proving the comparison principle.

\vspace{1mm}

\noindent {\bf Proof of Theorem \ref{Thmcompeps}. }

\vspace{1mm}

\noindent Let $u$ and $w$ as in Theorem \ref{Thmcompeps}, and (re)define $w$ on $[0,T]\times\partial\Sc_\eps$ by 
\beq \label{defwparsc}
w(t,z,\theta) &=& \liminf_{\tiny{\begin{array}{c}(t',z',\theta') \rightarrow (t,z,\theta)\\ (t',z',\theta') \in [0,T)\times\Sc_\eps\end{array}}} w(t',z',\theta'), 
\;\;\;\;\; (t,z,\theta) \in [0,T] \times \partial\Sc_\eps. 
\enq
In order to obtain the comparison result \reff{uleqw}, it suffices to prove that $\sup_{[0,T]\times\bar\Sc_\eps}(u-w)$ $\leq$ $0$, and we shall argue  
by contradiction by assuming  that 
\beq \label{contrary}
\sup_{[0,T]\times\bar\Sc_\eps} (u-w) & > & 0.
\enq

\noindent $\bullet$ {\it Step 1. Construction of a strict viscosity supersolution.}  

\vspace{1mm}

\noindent  Consider the  function defined  on $[0,T]\times\bar\Sc_\eps$ by
\beqs
\psi(t,z,\theta) &=& e^{\rho'(T-t)}  L_M(z)^{\gamma'} , \;\;\;\;\;   t \in [0,T], (z,\theta)=(x,y,p,\theta) \in \bar\Sc_\eps,
\enqs
where   $\rho'$ $>$ $0$, and $\gamma'$ $\in$ $(0,1)$ will be chosen later.  The function $\psi$ is  smooth $C^2$ on 
$[0,T)\times(\bar\Sc_\eps\setminus D_0)$, and by the same calculations as in \reff{derivLM}, we see that  by choosing 
$\rho'$ $>$ $\frac{\gamma'}{1-\gamma'}\frac{b^2}{2\sigma^2}$, then 
\beq \label{psiLc}
- \Dt{\psi} - \Dth{\psi} - \Lc \psi & > & 0 \;\;\;\; \mbox{ on } [0,T)\times(\bar\Sc_\eps\setminus D_{0}). 
\enq
Moreover, from \reff{HcLM}, we have  
\beq
(\psi - \Hc_\eps \psi )(t,z,\theta) &=& e^{\rho'(T-t)} \Big[L_M(z)^{\gamma'}  -  (\Hc_\eps L_M(z))^{\gamma'} \Big]  \; =: \; \Delta(t,z) 
\label{defDelta} \\
& > & 0   \;\;\;\;\;\;\;   \mbox{ on } [0,T] \times  \bar\Sc_\eps.  \nonumber 
\enq
For $m$ $\geq$ $1$, we denote by 
\beqs
\tilde u(t,z,\theta) \; = \; e^t u(t,z,\theta), & \mbox{ and } & \tilde w_m(t,z,\theta) \;=\;  e^t \big[ w(t,z,\theta) + \frac{1}{m} \psi(t,z,\theta)]. 
\enqs
From the viscosity subsolution property of $u$,  we immediately see that $\tilde u$ is a viscosity subsolution to 
\beq
\min\big[ \tilde u - \Dt{\tilde u} - \Dth{\tilde u} - \Lc \tilde u \; , \; \tilde u - \Hc_\eps\tilde u\big] & \leq & 0, \;\;\; \mbox{ on } \; [0,T)\times\bar\Sc_\eps 
\label{viscotildeu} \\
\min\big[ \tilde u -  \tilde U_{L_\eps} , \; \tilde u - \Hc_\eps\tilde u\big] & \leq & 0, \;\;\; \mbox{ on }  \; \{T\} \times\bar\Sc_\eps, \label{viscotildeuT}
\enq
where we set $\tilde U_{L_\eps}(z,\theta)$ $=$ $e^T U_{L_\eps}(z,\theta)$.  From the viscosity supersolution property of $w$, and the relations 
\reff{psiLc}-\reff{defDelta},  we also derive  that  $\tilde w_m$ is a viscosity supersolution to 
\beq
\tilde w_m  - \Dt{\tilde w_m} - \Dth{\tilde w_m} - \Lc \tilde w_m & \geq & 0 \;\;\;\;\;\;\;  \mbox{ on } \; [0,T)\times (\Sc_\eps\setminus D_{0}) 
\label{viscotildewm} \\
\tilde w_m - \Hc_\eps \tilde w_m & \geq & \frac{1}{m} \Delta \;\;\; \mbox{ on } \; [0,T] \times \Sc_\eps. \label{viscotildewmHc} \\
\tilde w_m - \tilde U_{L_\eps} & \geq & 0 \;\;\;\;\;\;\;  \mbox{ on } \; \{T\}\times\Sc_\eps. \label{viscotildewmT}
\enq
On the other hand,  from the growth condition on $u$ and $w$ in  $\Gc_{\gamma}([0,T]\times\bar\Sc _{\eps})$, and by choosing $\gamma' $ $\in$ 
$(\gamma,1)$, we have for all $(t,\theta)$ $\in$ $[0,T]^2$, 
\beqs
\lim_{|z| \rightarrow \infty} (u-w_m)(t,z,\theta) &=& - \infty. 
\enqs
Therefore, the usc function $\tilde u-\tilde w_m$ attains its supremum on $[0,T]\times\bar\Sc_\eps$, and from \reff{contrary}, there exists  
$m$ large enough, and $(\bar t,\bar z,\bar\theta)$  $\in$ $[0,T]\times\bar\Sc_\eps$ s.t. 
\beq \label{deftildeM}
\tilde M \; = \;  \sup_{[0,T]\times\bar\Sc_\eps} (\tilde u -\tilde w_m) &=&  (\tilde u-\tilde w_m)(\bar t,\bar z,\bar\theta) \; > \; 0.
\enq

\noindent $\bullet$ {\it Step 2.}  From the boundary condition \reff{uwboun}, we know that $(\bar z,\bar\theta)$ cannot lie in $D_0$, and 
we have then two possible cases: 

(i) $(\bar z,\bar \theta)$ $\in$ $\Sc_\eps\setminus D_0$

\vspace{1mm}

(ii)  $(\bar z,\bar \theta)$ $\in$ $\partial\Sc_\eps\setminus D_0$.

\vspace{1mm}

\noindent The case (i) where $(\bar z,\bar \theta)$ lies in  $\Sc_\eps$ 
is standard in the comparison principle for (nonconstained) viscosity solutions, and we focus here on the case 
(ii), which is specific  to  cons\-trained viscosity solutions.  From \reff{defwparsc},  there exists a sequence 
$(t_{n},z_{n},\theta_n)_{n\geq1}$ in $[0,T)\times\Sc _{\eps}$ such that 
\beqs
(t_{n},z_{n},\theta_n,\tilde w_{m}(t_{n},z_{n},\theta_n)) & \longrightarrow &  (\bar t,\bar z,\bar\theta,\tilde w_{m}(\bar t,\bar z,\bar\theta)) 
\qquad \mbox{ as } ~ n \rightarrow\infty.
\enqs
We then set   $\delta_{n}$ $=$ $|z_{n}- \bar z|+|\theta_n-\bar\theta|$, and consider the function 
$\Phi_{n}$ defined on $[0,T]\times(\bar\Sc _{\eps})^2$ by:
\beqs
\Phi_{n}(t,z,\theta,z',\theta') & = & \tilde u(t,z,\theta)- \tilde w_{m}(t,z',\theta')- \varphi_{n}(t,z,\theta,z',\theta') \\
\varphi_{n}(t,z,\theta,z',\theta') & = & |t- \bar t|^2 + |z- \bar z|^4  + |\theta - \bar\theta|^4  \\
& & \; + \;    \frac{|z-z'|^2 + |\theta-\theta'|^2}{2\delta_{n}} + \Big( \frac{d(z',\theta')}{d(z_{n},\theta_n)}-1\Big)^4. 
\enqs
Here, $d(z,\theta)$ denotes the distance from $(z,\theta)$ to $\partial \Sc _{\eps}$. Since $(\bar z,\bar\theta)$ $\notin$ $D_0$, 
there exists an open neighborhood $\bar\Vc$  of $(\bar z,\bar\theta)$ satisfying $\bar\Vc\cap D_{0} = \emptyset$, such that the function $d(.)$ is twice continuously differentiable with bounded derivatives. This is well known (see e.g. \cite{giltru77}) when $(\bar z,\bar \theta)$ lies in  the smooth  parts  of the boundary  
$\partial\Sc_\eps\setminus( D_{1,\eps}\cup D_{2,\eps})$.  This is also true for $(\bar z,\bar\theta)$ $\in$ $D_{k,\eps}$ for $k\in\{1,2\}$.  Indeed, at these corner lines, the inner normal vectors form an acute angle (positive scalar product), and thus one can extend  from $(\bar z,\bar\theta)$ the boundary to a smooth boundary so that the distance $d$ is equal, locally on the neighborhood, to the distance to this  smooth boundary.    
From the growth conditions on $u$ and $w$ in  $\Gc_{\gamma}([0,T]\times\bar\Sc _{\eps})$, there exists a sequence 
$(\hat t_n,\hat z_n,\hat\theta_n,\hat z_n',\hat\theta_n')$ attaining the maximum of the usc $\Phi_n$ on $[0,T]\times(\bar\Sc_\eps)^2$. By standard arguments (see e.g. \cite{bar94} or \cite{lyvmnipha07}), we have
\beq
(\hat t_n,\hat z_n,\hat\theta_n,\hat z_n',\hat\theta_n') & \longrightarrow & (\bar t,\bar z,\bar\theta,\bar z,\bar\theta) \label{limtn} \\
\frac{|\hat z_n-\hat z_n'|^2 + |\hat \theta_n-\hat \theta_n'|^2}{2\delta_{n}} 
+ \Big( \frac{d(\hat z_n',\hat \theta_n')}{d(z_{n},\theta_n)}-1\Big)^4 & \longrightarrow & 0 \label{limdn}\\
\tilde u(\hat t_n,\hat z_n,\hat \theta_n) - \tilde w_m(\hat t_n,\hat z_n',\hat\theta_n') & \longrightarrow & (\tilde u -\tilde w_m)(\bar t,\bar z,\bar\theta). 
\label{limun}
\enq
The convergence in \reff{limdn} shows in particular that for $n$ large enough, $d(\hat z_n',\hat \theta_n')$ $\geq$ $d(z_n,\theta_n)/2$ $>$ $0$, and 
so $(\hat z_n',\hat\theta_n')$ $\in$ $\Sc_\eps$.  From the convergence in \reff{limtn}, we may also assume that for $n$ large enough, 
$(\hat z_n,\hat\theta_n)$, $(\hat z_n',\hat\theta_n')$ lie in the neighborhood $\bar\Vc$ of $(\bar z,\bar\theta)$ so that the derivatives 
upon order $2$ of $d(.)$ at  $(\hat z_n,\hat\theta_n)$ and  $(\hat z_n',\hat\theta_n')$ exist and are bounded.

\vspace{1mm}

\noindent $\bullet$ {\it Step 3.}  We show that for $n$ large enough, 
\beq \label{vstrict}
\tilde u(\hat t_n,\hat z_n,\hat \theta_n) -  \Hc_\eps \tilde u(\hat t_n,\hat z_n) & > & 0. 
\enq 
Otherwise, up to a subsequence, we would have for all $n$:  
\beqs
\tilde u(\hat t_n,\hat z_n,\hat \theta_n) -  \Hc_\eps \tilde u(\hat t_n,\hat z_n) & \leq & 0.  
\enqs
By sending $n$ to infinity, and from the upper-semicontinuity of $\Hc_\eps\tilde u$, we get with \reff{limtn}: 
$-\infty$ $<$ $\tilde u(\bar t,\bar z,\bar\theta)$ $\leq$ $\Hc_\eps \tilde u(\bar t,\bar z,\bar\theta)$, which shows in particular that 
$\Cc_\eps(\bar z,\bar\theta)$ is not empty. 
Moreover, by the viscosity supersolution property \reff{viscotildewmHc}, we have 
\beqs
 \tilde w_m(\hat t_n,\hat z_n',\hat\theta_n') - \Hc_\eps \tilde w_m(\hat t_n,\hat z_n',\hat\theta_n') & \geq & 
 \frac{1}{m} \Delta(\hat t_n,\hat z_n',\hat\theta_n'). 
 \enqs
By substracting the two previous inequalities, we would get
\beqs
\tilde u(\hat t_n,\hat z_n,\hat \theta_n) -   \tilde w_m(\hat t_n,\hat z_n',\hat\theta_n') & \leq &  \Hc_\eps \tilde u(\hat t_n,\hat z_n) 
 - \Hc_\eps \tilde w_m(\hat t_n,\hat z_n',\hat\theta_n') -  \frac{1}{m} \Delta(\hat t_n,\hat z_n',\hat\theta_n'). 
\enqs
By sending $n$ to infinity, and from the upper-semicontinuity of $\Hc_\eps\tilde u$,  the lower-semicontinuity of $\Hc_\eps \tilde w_m$ and 
$\Delta$,  this  yields with \reff{limtn}, \reff{limun}
\beqs
(\tilde u-\tilde w_m)(\bar t,\bar z,\bar\theta) & \leq & 
\Hc_\eps \tilde u(\bar t,\bar z,\bar\theta) -  \Hc_\eps \tilde w_m(\bar t,\bar z,\bar\theta) -  \frac{1}{m} \Delta(\bar t,\bar z,\bar\theta). 
\enqs
Now, by compactness of $\Cc_\eps(\bar z,\bar\theta)$ $\neq$ $\emptyset$, there exists $\bar e$ $\in$ $\Cc_\eps(\bar z,\bar\theta)$ such that 
$\Hc_\eps\tilde u(\bar t,\bar z,\bar\theta)$ $=$ $\tilde u(t,\Gamma_\eps(\bar z,\bar\theta,\bar e),0)$ and so
\beqs
\tilde M \; = \;  (\tilde u-\tilde w_m)(\bar t,\bar z,\bar\theta)  & \leq &  
\tilde u(\bar t,\Gamma_\eps(\bar z,\bar\theta,\bar e),0) - \tilde w_m(\bar t,\Gamma_\eps(\bar z,\bar\theta,\bar e),0)  -  
\frac{1}{m} \Delta(\bar t,\bar z,\bar\theta) \\
& \leq & \tilde M  - \frac{1}{m} \Delta(\bar t,\bar z,\bar\theta),
\enqs
a contradiction. 

\vspace{1mm}

\noindent $\bullet$ {\it Step 4.} We check that, up to a subsequence, $\hat t_n$ $<$ $T$ for all $n$.   On the contrary, $\hat t_n$ $=$ $\bar t$ $=$ $T$ 
for $n$ large enough, and we would get   from \reff{vstrict} and the viscosity subsolution property \reff{viscotildeuT}: 
\beqs
\tilde u(T,\hat z_n,\hat\theta_n) & \leq & \tilde U_{L_\eps}(\hat z_n,\hat\theta_n). 
\enqs
Moreover,  by \reff{viscotildewmT}, we have $\tilde w_m(T,\hat z_n',\hat\theta_n')$ $\geq$  $\tilde U_{L_\eps}(\hat z_n',\hat\theta_n')$, which combined with the former  inequality, implies 
\beqs
\tilde u(T,\hat z_n,\hat\theta_n) -  \tilde w_m(T,\hat z_n',\hat\theta_n') & \leq & 
\tilde U_{L_\eps}(\hat z_n,\hat\theta_n) -  \tilde U_{L_\eps}(\hat z_n',\hat\theta_n'). 
\enqs
By sending $n$ to infinity,  this yields  with \reff{limtn}, \reff{limun} and continuity of $\tilde U_{L_\eps}$: 
$\tilde M$ $=$ $(\tilde u-\tilde w_m)(\bar t,\bar z,\bar\theta) $ $\leq$ $0$, a contradiction with \reff{deftildeM}. 

\vspace{1mm}

\noindent $\bullet$ {\it Step 5.}  We use the viscosity subsolution property \reff{viscotildeu} of $\tilde u$ at $(\hat t_n,\hat z_n,\hat \theta_n)$ $\in$ 
$[0,T)\times\bar\Sc_\eps$, which is written by   \reff{vstrict} as
\beq \label{ishiitildeu}
(\tilde u - \Dt{\tilde u} - \Dth{\tilde u} - \Lc \tilde u)(\hat t_n,\hat z_n,\hat\theta_n) & \leq & 0. 
\enq
The above inequality is understood in the viscosity sense, and applied with the test function 
$(t,z,\theta)$ $\rightarrow$ $\varphi_n(t,z,\theta,\hat z_n',\hat\theta_n')$, which is $C^2$ in the neighborhood $[0,T]\times\bar\Vc$ of 
$(\hat t_n,\hat z_n,\hat\theta_n)$.   We also write  the viscosity supersolution property \reff{viscotildewm} of $\tilde w_m$ at  
$(\hat t_n,\hat z_n',\hat\theta_n')$ $\in$ $[0,T)\times(\Sc_\eps\backslash
 D_{0})$:
\beq \label{ishiitildewm}
(\tilde w_m - \Dt{\tilde w_m} - \Dth{\tilde w_m} - \Lc \tilde w_m)(\hat t_n,\hat z_n',\hat\theta_n') & \geq & 0. 
\enq 
The above inequality is again understood in the viscosity sense, and applied with the test function 
$(t,z',\theta')$ $\rightarrow$ $-\varphi_n(t,\hat z_n,\hat \theta_n,z',\theta')$, which is $C^2$ in the neighborhood $[0,T]\times\bar\Vc$ of 
$(\hat t_n,\hat z_n',\hat\theta_n')$.   The conclusion is achieved by  arguments  similar to \cite{lyvmnipha07}:  
we invoke Ishii's Lemma, substract the two inequalities \reff{ishiitildeu}-\reff{ishiitildewm},  and  finally  get the required contradiction $\tilde M$ 
$\leq$ $0$ by sending $n$ to infinity with \reff{limtn}-\reff{limdn}-\reff{limun}. 
\ep

\section{An approximating problem with utility penalization}

\setcounter{equation}{0} \setcounter{Assumption}{0}
\setcounter{Theorem}{0} \setcounter{Proposition}{0}
\setcounter{Corollary}{0} \setcounter{Lemma}{0}
\setcounter{Definition}{0} \setcounter{Remark}{0}

We consider in this section another perturbation  of our initial  optimization problem by  adding a cost $\eps$ to the  utility  at each trading. 
We then define the value function $\bar v_{\eps}$ on $[0,T]\times\bar\Sc$ by 
\beq\label{PenU}
\bar v_{\eps}(t,z,\theta) & = & \sup_{\alpha\in\Ac_{\ell}^b(t,z,\theta)}\E\Big[U_{L}\big(Z_{T},\Theta_{T}\big)-\eps N_{T}(\alpha)\Big], 
\;\;\; (t,z,\theta) \in [0,T]\times\bar\Sc. 
\enq

The convergence of this approximation is immediate.

\begin{Proposition} \label{propconv}
The sequence $(\bar v_{\eps})_{\eps}$ is nondecreasing and converges pointwise on $[0,T]\times\bar \Sc$ towards $v$ as $\eps$ goes to zero.
\end{Proposition}
\textbf{Proof.} 
It is clear that the sequence $(\bar v_{\eps})_{\eps}$ is nondecreasing and that 
$\bar v_\eps\leq v$ on $[0,T]\times\bar \Sc$ for any $\eps>0$. Let us prove that $\lim_{\eps\searrow0}\bar v_\eps$ $=$ $v$. 
Fix $n\in\N^*$ and $(t,z,\theta)\in[0,T]\times\bar \Sc$ and consider some $\alpha^{(n)}\in\Ac^b_{\ell}(t,z,\theta)$ such that 
\beqs
\E\Big[U_{L}\big(Z^{(n)}_{T},\Theta^{(n)}_{T}\big)\Big] & \geq & v(t,z,\theta)-\frac{1}{n}, 
\enqs
where $(Z^{(n)},\Theta^{(n)})$ is the associated controlled process. From the monotone convergence theorem, we then get
\beqs
\lim_{\eps\searrow0} \bar v_\eps(t,z,\theta) & \geq & \E\Big[U_{L}\Big(Z^{(n)}_{T},\Theta^{(n)}_{T}\Big)\Big]  ~\geq~ v(t,z,\theta)-\frac{1}{n}.
\enqs 
By the arbitrariness of $n\in \N^*$, we  conclude that   $\lim_\eps \bar v_\eps$ $\geq$ $v$, which ends the proof since we already have 
$\bar v_\eps$ $\leq$ $v$. 
\ep

\vspace{2mm}

The nonlocal impulse  operator $\bar\Hc^\eps$ associated to \reff{PenU} is given by
\beqs
\bar \Hc_\eps\varphi(t,z,\theta) & = & \Hc\varphi(t,z,\theta)-\eps,
\enqs
and we consider the corresponding dynamic programming equation:
\beq
 \min\Big[ -\Dt{w} - \Dth{w} - \Lc w~,~w -\bar \Hc_\eps w \Big] & = & 0, \;\;\;\; \mbox{ in } [0,T)\times\bar\Sc, \label{QVIUeps} \\
 \min\big[ w -  U_{L}~,~w - \bar \Hc_\eps w \big] &=& 0,  \;\;\;\; \mbox{ in }  \{T\} \times\bar\Sc. \label{EDPBUeps}
\enq

\vspace{2mm}

By similar arguments as in Section 5, we can show that $\bar v_\eps$ is a constrained viscosity solution to  
\reff{QVIUeps}-\reff{EDPBUeps}, and  the following comparison principle holds:

\noindent Suppose $u$  $\in$ $\Gc_{\gamma}([0,T]\times\bar\Sc )$ is a usc viscosity subsolution to  \reff{QVIUeps}-\reff{EDPBUeps} on 
$[0,T]\times\bar\Sc$, and $w$ $\in$ $\Gc_{\gamma}([0,T]\times\bar\Sc )$ is a lsc  viscosity supersolution to  \reff{QVIUeps}-\reff{EDPBUeps} on 
$[0,T]\times\Sc$, such that 
\beqs
u(t,z,\theta) & \leq & \liminf_{\tiny{\begin{array}{c}(t',z',\theta')\rightarrow(t,z,\theta)\\ (t',z',\theta')\in[0,T)\times\Sc \end{array}}}w(t',z',\theta'), 
\quad\forall(t,z,\theta)\in [0,T] \times D_0. 
\enqs
Then, 
\beq \label{uleqw2}
u & \leq & w \quad\mbox{ on }\quad [0,T]\times \Sc. 
\enq
The proof follows the same lines of arguments as in the proof of Theorem \ref{Thmcompeps} 
(the function $\psi$ is still a strict viscosity supersolution to  \reff{QVIUeps}-\reff{EDPBUeps} on $[0,T]\times\bar \Sc$), and so we omit it.

As a consequence, we obtain  a  PDE characterization of the value function $v$.

\begin{Proposition}
The value function $v$ is the minimal constrained viscosity solution in  $\Gc_{\gamma}([0,T]\times\bar\Sc )$ 
to \reff{QVIv}-\reff{condterm}, satisfying the boundary condition 
\beq \label{bounv}
\lim_{(t',z',\theta')\rightarrow (t,z,\theta)} v(t',z',\theta') &=& v(t,z,\theta) \; = \;  U(0), \;\;\;  \forall (t,z,\theta) \in [0,T]\times D_0.  
\enq
\end{Proposition}
\textbf{Proof.} Let $V\in \Gc_{\gamma}([0,T]\times\bar\Sc )$ be a viscosity solution in  $\Gc_{\gamma}([0,T]\times\bar\Sc )$ 
to \reff{QVIv}-\reff{condterm}, satisfying the boundary condition \reff{bounv}.  Since $\Hc$ $\geq$ $\bar\Hc_\eps$, it is clear that 
$V_*$ is a viscosity supersolution  to \reff{QVIUeps}-\reff{EDPBUeps}.   Moreover,  since 
$\lim_{(t',z',\theta')\rightarrow (t,z,\theta)} V_*(t',z',\theta')$ $=$ $U(0)$ $=$ $v(t,z,\theta)$ $\geq$ $\bar v_\eps^*(t,z,\theta)$ for 
$(t,z,\theta)$ $\in$ $[0,T]\times D_0$, we deduce from the comparison principle \reff{uleqw2} that 
$V$ $\geq$ $V_{*}$ $\geq$ $\bar v_\eps^*$ $\geq$ $\bar v_{\eps}$ on $[0,T]\times \Sc$. 
By sending $\eps$ to 0, and from the convergence result in Proposition \ref{propconv}, we obtain: $V$ $\geq$ $v$, which proves the required result.  
\ep

\section*{Appendix:  constrained viscosity solutions to parabolic  QVIs}

\renewcommand{\theLemma}{A.\arabic{Lemma}}
\renewcommand{\theDefinition}{A.\arabic{Definition}}
\renewcommand{\theequation}{A.\arabic{equation}}
\renewcommand{\thesection}{A.\arabic{section}}

\setcounter{section}{0}
\setcounter{equation}{0}
\setcounter{Lemma}{0}
\setcounter{Definition}{0}

We consider a parabolic quasi-variational inequality in the form: 
\beq \label{QVIgen}
\min \Big[ - \Dt{v}  + F(t,x,v,D_x v,D_{x}^2 v) \;,\; v - \Hc v \Big] &=& 0, \;\;\; \mbox{ in } [0,T)\times\bar\Oc,
\enq
together with a terminal condition
\beq \label{termgen}
\min \big[ v - g \; , \; v - \Hc v \big] &=& 0,  \;\;\; \mbox{ in }  \{T\}\times\bar\Oc.
\enq
Here, $\Oc$ $\subset$ $\R^d$ is an open domain, $F$ is a continuous function on 
$[0,T]\times\R^d\times\R\times\R^d\times\Sc^d$ $(\Sc^d$ is the set of positive semidefinite  symmetric  matrices in $\R^{d\times d}$), 
nonincreasing in its last argument, $g$ is a continuous function on $\bar\Oc$, 
and $\Hc$ is a nonlocal operator defined on the set of locally bounded functions on $[0,T]\times\bar\Oc$ by: 
\beqs
\Hc  v(t,x) &=& \sup_{e \in \Cc(t,x)} \big[  v(t, \Gamma(t,x,e)) + c(t,x,e) \big]. 
\enqs 
$\Cc(t,x)$ is a compact set of a metric space $E$, eventually empty for some values of $(t,x)$, in which case we set $\Hc v(t,x)$ $=$ $-\emptyset$, 
and is  continuous for the Hausdorff metric, i.e. if $(t_n,x_n)$ converges to $(t,x)$ in 
$[0,T]\times\bar\Oc$, and $(e_n)$ is a sequence in $\Cc(t_{n},x_n)$ converging to $e$, then $e$ $\in$ $\Cc(t,x)$. 
The  functions $\Gamma$ and $c$ are continuous, and such that  $\Gamma(t,x,e)$ $\in$ $\bar\Oc$ for all $e$ $\in$ $\Cc(t,x,e)$.

Given a locally bounded function $u$ on $[0,T]\times\bar\Oc$, we define  its lower-semicontinuous (lsc in short) envelope $u_*$ and 
upper-semicontinuous (usc) envelope $u^*$ on $[0,T]\times\bar\Sc$ by: 
\beqs
u_*(t,x) \; = \;  \liminf_{\tiny{\begin{array}{c}(t',x')\rightarrow (t,x)\\ (t',x')\in[0,T)\times\Oc\end{array}}} u(t',x'), & &  u^*(t,x) \; = \;  \limsup_{\tiny{\begin{array}{c}(t',x')\rightarrow (t,x)\\ (t',x')\in[0,T)\times\Oc\end{array}}} u(t',x'). 
\enqs
One can  check (see e.g. Lemma 5.1 in \cite{lyvmnipha07}) that the operator $\Hc$ preserves lower and upper-semicontinuity: 

\vspace{1mm}

(i) $\Hc u_*$ is lsc, and  $\Hc u_*$ $\leq$ $(\Hc u)_*$, \hspace{7mm} (ii) $\Hc u^*$ is usc, and $(\Hc u)^*$ $\leq$ $\Hc u^*$. 
 
\vspace{1mm}

We now give the definition of constrained viscosity solutions  to \reff{QVIgen}-\reff{termgen}. This notion, which extends the definition of 
viscosity solutions of Crandall, Ishii and Lions (see \cite{craishlio92}), was introduced in \cite{son86} for first-order equations for taking into account boundary conditions arising in state constraints, and used in \cite{zar88} for stochastic control problems in optimal investment.

\begin{Definition}
A locally bounded function $v$ on $[0,T]\times\bar\Oc$ is a constrained viscosity solution to \reff{QVIgen}-\reff{termgen} if the two following properties hold:

\vspace{1mm}

\noindent (i) Viscosity supersolution property on $[0,T]\times\Oc$: for all $(\bar t,\bar x)$ $\in$ $[0,T]\times\Oc$, and $\varphi$ $\in$ 
$C^{1,2}([0,T]\times\Oc)$ with $0$ $=$ $(v_*-\varphi)(\bar t,\bar x)$ $=$ $\min(v_*-\varphi)$, we have
\beqs
\min \Big[ - \Dt{\varphi}(\bar t,\bar x)   
+ F(\bar t,\bar x,\varphi_*(\bar t,\bar x),D_x \varphi(\bar t,\bar x),D_{x}^2 \varphi(\bar t,\bar x)) \;  , & &  \\
  \;\;\; v_*(\bar t,\bar x) - \Hc v_*(\bar t,\bar x) \Big] & \geq &  0, \;\; (\bar t,\bar x) \in  [0,T)\times\Oc, \\
\min \big[ v_*(\bar t,\bar x)  - g(\bar x) \; , \;  v_*(\bar t,\bar x) - \Hc v_*(\bar t,\bar x) \Big] & \geq &  0, \;\;\; (\bar t,\bar x) \in \{T\}\times\Oc.   
\enqs

\noindent (ii) Viscosity subsolution property on $[0,T]\times\bar\Oc$: for all $(\bar t,\bar x)$ $\in$ $[0,T]\times\bar\Oc$, and $\varphi$ $\in$ 
$C^{1,2}([0,T]\times\bar\Oc)$ with $0$ $=$ $(v^*-\varphi)(\bar t,\bar x)$ $=$ $\max(v^*-\varphi)$, we have
\beqs
\min \Big[ - \Dt{\varphi}(\bar t,\bar x)   
+ F(\bar t,\bar x,\varphi_*(\bar t,\bar x),D_x \varphi(\bar t,\bar x),D_{x}^2 \varphi(\bar t,\bar x)) \;  , & &  \\
  \;\;\; v^*(\bar t,\bar x) - \Hc v^*(\bar t,\bar x) \Big] & \leq &  0, \;\; (\bar t,\bar x) \in  [0,T)\times\bar\Oc, \\
\min \big[ v^*(\bar t,\bar x)  - g(\bar x) \; , \;  v^*(\bar t,\bar x) - \Hc v^*(\bar t,\bar x) \Big] & \leq &  0, \;\;\; (\bar t,\bar x) \in \{T\}\times\bar\Oc.   
\enqs
\end{Definition}

\vspace{7mm}

\begin{small}

\end{small}

\end{document}